\documentclass[11pt, reqno]{amsart} 

\usepackage{fancyhdr}
\usepackage{amsmath}
\usepackage{amsthm}
\usepackage{amssymb}
\usepackage[foot]{amsaddr}
\usepackage{enumerate}
\usepackage{tikz}
\usetikzlibrary{calc}
\usepackage{mathtools}
\usepackage[margin=1in]{geometry}

\usepackage[unicode]{hyperref}
\usepackage{xcolor}
\hypersetup{
	colorlinks,
    linkcolor={blue!50!black},
    citecolor={blue!50!black},
    urlcolor={blue!80!black}}

\title{Continuous limits of generalized pentagram maps}

\author{Danny Nackan$^*$}
\address{$^*$Department of Mathematics, Yale University, New Haven, CT 06511, USA}
\email{$^*$danny.nackan@yale.edu}

\author{Romain Speciel$^\dagger$}
\address{$^\dagger$Department of Mathematics, University of Toronto, Toronto, ON M5S 2E4, Canada}
\email{$^\dagger$romain.speciel@mail.utoronto.ca}

\makeatletter
\let\mytitle\@title
\makeatother

\newtheorem{theorem}{Theorem}[section]
\newtheorem{lemma}[theorem]{Lemma}
\newtheorem{cor}[theorem]{Corollary}
\newtheorem{prop}[theorem]{Proposition}

\theoremstyle{definition}
\newtheorem{definition}[theorem]{Definition}
\newtheorem{example}[theorem]{Example}
\newtheorem{remark}[theorem]{Remark}

\newtheorem{examplecont}{Example}

\newtheorem{theoremalph}{Theorem}

\newcommand{\R}{\mathbb{R}}
\newcommand{\RP}{\mathbb{RP}}
\newcommand{\Z}{\mathbb{Z}}

\newcommand{\e}{\epsilon}
\renewcommand{\P}{\mathcal{P}}
\renewcommand{\O}{\mathcal{O}}

\DeclareMathOperator{\PGL}{PGL}

\DeclareMathOperator{\DO}{DO}

\DeclareMathOperator{\id}{id}

\renewcommand{\epsilon}{\varepsilon}

\newcommand{\mat}[1]{\begin{pmatrix} #1 \end{pmatrix}}	
\renewcommand{\l}{\left}
\renewcommand{\r}{\right}								
\newcommand{\tn}[1]{\textnormal{#1}}					
\renewcommand{\ge}{\Gamma_\e}

\makeatletter
\renewcommand{\tocsection}[3]{%
  \indentlabel{\@ifnotempty{#2}{\bfseries\ignorespaces#1 #2\quad}}\bfseries#3}
\renewcommand{\tocsubsection}[3]{%
  \indentlabel{\@ifnotempty{#2}{\ignorespaces#1 #2\quad}}#3}

\newcommand\@dotsep{4.5}
\def\@tocline#1#2#3#4#5#6#7{\relax
  \ifnum #1>\c@tocdepth 
  \else
    \par \addpenalty\@secpenalty\addvspace{#2}%
    \begingroup \hyphenpenalty\@M
    \@ifempty{#4}{%
      \@tempdima\csname r@tocindent\number#1\endcsname\relax
    }{%
      \@tempdima#4\relax
    }%
    \parindent\z@ \leftskip#3\relax \advance\leftskip\@tempdima\relax
    \rightskip\@pnumwidth plus1em \parfillskip-\@pnumwidth
    #5\leavevmode\hskip-\@tempdima{#6}\nobreak
    \leaders\hbox{$\m@th\mkern \@dotsep mu\hbox{.}\mkern \@dotsep mu$}\hfill
    \nobreak
    \hbox to\@pnumwidth{\@tocpagenum{\ifnum#1=1\bfseries\fi#7}}\par
    \nobreak
    \endgroup
  \fi}
\AtBeginDocument{%
\expandafter\renewcommand\csname r@tocindent0\endcsname{0pt}
}
\def\l@subsection{\@tocline{2}{0pt}{2.5pc}{5pc}{}}
\makeatother

\begin{document}


\begin{abstract}
    We provide a rigorous  treatment of continuous limits for various generalizations of the pentagram map on polygons in $\RP^d$ 
    by means of quantum calculus. Describing this limit in detail for the case of the short-diagonal pentagram map, we verify that this construction yields the $(2,d+1)$-KdV equation, and moreover, the Lax form of the pentagram map in the limit is proved to become the Lax representation 
    of the corresponding KdV system.
    More generally, we introduce the $\chi$-pentagram map, a geometric construction defining curve evolutions by directly taking intersections of subspaces through specified points. We show that its different configurations yield certain other KdV equations and provide an argument towards disproving the conjecture that any KdV-type equation can be discretized through pentagram-type maps.
\end{abstract}

\date{} 
\maketitle
\thispagestyle{empty}

\tableofcontents


\section{Introduction}

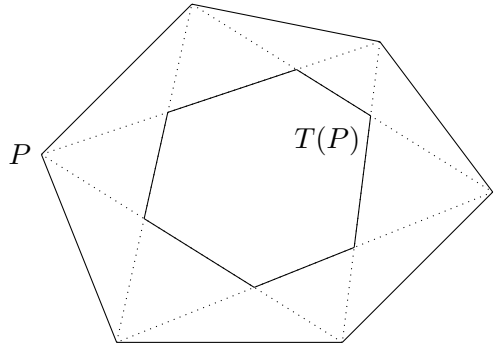
\begin{figure}[ht]
\begin{center}
	\begin{equation*}
		\begin{tikzpicture} 
		\coordinate (v1) at (-2,1);
		\coordinate (v2) at (0,3);
		\coordinate (v3) at (2.5,2.5);
		\coordinate (v4) at (4,0.5);
		\coordinate (v5) at (2,-1.5);
		\coordinate (v6) at (-1, -1.5);
		\draw (v1) -- (v2);
		\draw (v2) -- (v3);
		\draw (v3) -- (v4);
		\draw (v4) -- (v5);
		\draw (v5) -- (v6);
		\draw (v6) -- (v1);
		\draw[dotted] (v1) -- (v3);
		\draw[dotted] (v2) -- (v4);
		\draw[dotted] (v3) -- (v5);
		\draw[dotted] (v4) -- (v6);
		\draw[dotted] (v5) -- (v1);
		\draw[dotted] (v6) -- (v2);
		\coordinate (u1) at (intersection of {v1--v3} and {v2--v4});
		\coordinate (u2) at (intersection of {v2--v4} and {v3--v5});
		\coordinate (u3) at (intersection of {v3--v5} and {v4--v6});
		\coordinate (u4) at (intersection of {v4--v6} and {v5--v1});
		\coordinate (u5) at (intersection of {v5--v1} and {v6--v2});
		\coordinate (u6) at (intersection of {v6--v2} and {v1--v3});
		\draw (u1) -- (u2);
		\draw (u2) -- (u3);
		\draw (u3) -- (u4);
		\draw (u4) -- (u5);
		\draw (u5) -- (u6);
		\draw (u6) -- (u1);
		\node[left] at (v1) {$P$};
		\node[below left] at (u2) {$T(P)$};
		\end{tikzpicture}
	\end{equation*}
\end{center}
\caption{The pentagram map ($n=6$)}
\label{fig:pentagram}
\end{figure}

The pentagram map $T$ was introduced by R.~Schwartz in~\cite{Schwartz:1992}: it sends a plane convex $n$-gon $P$ to the new $n$-gon $T(P)$ whose vertices are formed by the intersections of the shortest diagonals of~$P$ (see Figure~\ref{fig:pentagram}). 
This definition can naturally be extended to the more general space of twisted $n$-gons in $\RP^2$ modulo projective equivalence. In this context, Ovsienko, Schwartz, and Tabachnikov proved in~\cite{Ovsienko:2010} that the pentagram map is a discrete integrable system, and that its continuous limit is the classical Boussinesq equation. 
The \emph{short-diagonal} map introduced by Khesin and Soloviev in~\cite{Khesin:2013} further generalizes the pentagram map to an arbitrary dimension~$d$ by taking intersections of $d$ hyperplanes through vertices in $\RP^d$. 
It was shown in~\cite{Khesin:2013} that the short-diagonal map is also integrable, and that its continuous limit is the $(2,d+1)$-equation of the KdV hierarchy. 

There are many other generalizations of the pentagram map; for instance, 
Khesin and Soloviev investigate different configurations of hyperplanes in~\cite{Khesin:2016},  Mar{\'\i}~Beffa considers intersections of certain other types of subspaces in~\cite{Mari-Beffa:2013}, while Izosimov relates such maps to Poisson-Lie groups of pseudo-difference operators in \cite{Izosimov:2019}.
In this paper we focus on continuous limits of pentagram-type maps.

To this end, we introduce the \emph{$\chi$-pentagram map} $T_\epsilon^\chi$, which defines an evolution of a curve~$\gamma$ by taking the intersection of subspaces of the form $P_\e^i(x)=\big(\gamma(x+p_{i,0}\e),\dots, \gamma(x+p_{i,q_i}\e)\big)$, 
with parameters $\chi=\{p_{i,j}\}$.
This construction provides a continuous analogue of all of the aforementioned generalizations of the pentagram map, and it is in this general setting that we study their evolution as $\epsilon \to 0$.

Below, we begin with a detailed treatment of the continuous limit for the short-diagonal pentagram map before discussing the general case. 
In the standard constructions of the continuous limit used so far (and summarized in Section~\ref{sec:background}), twisted polygons were replaced by nondegenerate smooth curves with monodromy, and the analogue of the pentagram map was set to be the evolution of such a curve in the direction of the envelope curve of a family of hyperplanes, cf.~\cite{Ovsienko:2010, Khesin:2013, Izosimov:2019}.
We provide an alternative construction: rather than considering the envelope curve, we investigate the evolution described by directly applying the pentagram map to a discretized curve $\gamma$.
(This construction can be described by a corresponding configuration of the $\chi$-pentagram map.)
In this case, one can use quantum calculus (see~\cite{Kac:2002}) to rigorously study the limit, as follows:
\begin{theoremalph}[=\,Proposition~\ref{prop:limit-of-difference-eq} and Theorem~\ref{thm:directctslimit}]
    The differential equations defining nondegenerate curves with monodromy are the quantum calculus limits of the difference equations defining twisted $n$-gons, as $n \to \infty$. 
    Furthermore, in the limit, the dynamics of the $n$-gons under the pentagram map become the dynamics of the curves, as described by the $(2, d+1)$-KdV equation.
\end{theoremalph}

We provide a similar analysis of the Lax form of the short-diagonal map.
Notably, our construction of the continuous limit allows one to directly compute the limiting dynamics of the involved matrices, even though explicit formulas for the Lax representation are not fully known for $d > 3$.

\begin{theoremalph}[=\,Theorem~\ref{thm:lax_limit}]
    In the continuous limit as $n \to \infty$, the Lax matrices associated to twisted $n$-gons tend to differential operators associated to curves.
    Furthermore, in the limit, the dynamics of the discrete Lax matrices become the dynamics of the differential operators, as described by the $(2, d+1)$-KdV zero-curvature equation.
\end{theoremalph}

In the more general setting of the $\chi$-pentagram map, we show that a broad class of subspace configurations yield the $(2,d+1)$-KdV equation:

\begin{theoremalph}[=\,Theorem~\ref{thm:chi_pent_cont_limit}]
    Under suitable parametrization conditions (corresponding to a geometric centralization of the evolution), the continuous limit of the $\chi$-pentagram map corresponds to the $(2,d+1)$-KdV equation.
\end{theoremalph}
In particular, this allows one to compute the continuous limit of the dual dented pentagram map, defined in \cite{Khesin:2016} (see Proposition~\ref{prop:dual_dented_limit} below).

Finally, in Section~\ref{sec:kdv_realization} we consider the question of realizing other KdV equations as continuous limits, building directly from the results established in~\cite{Mari-Beffa:2013}.
We construct new instances of the $\chi$-pentagram map which yield the $(3,4)$-KdV equation as their continuous limit, using the methods developed in Section~\ref{sec:ctslimitcalc}. 
However, we conclude by providing a heuristic argument (=\,Proposition~\ref{prop:degrees_of_freedom_count}) which strongly suggests that not all KdV equations can be obtained as continuous limits of pentagram-type maps, contrary to a conjecture in~\cite{Mari-Beffa:2013}.

\vspace{1em}

\noindent\textbf{Acknowledgments.} This work was partially supported by the Natural Sciences and Engineering Research Council of Canada [USRA to R.S.]. Both authors wish to sincerely thank Prof.~Boris Khesin for the many helpful discussions and suggestions in preparing this paper. We also thank the anonymous referees for their comments.


\section{Background and notation} \label{sec:background}

\subsection{Short-diagonal pentagram maps}\label{sub:short_diag_background}

In this section we recall the description of the short-diagonal map, following~\cite{Ovsienko:2010} in the 2D case and 
\cite{Khesin:2013} in an arbitrary dimension $d$.

\begin{definition}
	A \emph{twisted $n$-gon} in $\RP^d$ is a map $\phi: \Z \to \RP^d$ such that $\phi(k+n) = M \cdot \phi(k)$ for every $k \in \Z$ and some $M \in \PGL(d+1, \R)$, called the \emph{monodromy} of $\phi$. 
\end{definition}
Note that when $M = \id$, we recover the usual notion of a closed $n$-gon with vertices $v_i\coloneqq \phi(i)$ in $\RP^d$, where $i=1, \dots, n$. Two twisted $n$-gons $\phi_1, \phi_2$ are said to be projectively equivalent if $\phi_2 = g \circ \phi_1$ for some $g \in \PGL(d+1, \R)$.
We will only consider those twisted $n$-gons such that any $d+1$ consecutive vertices are in general position. Let $\P_n$ denote the moduli space of such twisted $n$-gons, considered up to projective equivalence.

Given $\phi \in \P_n$, define for each $i \in \Z$ the \emph{short-diagonal hyperplane} $P_i$ passing through $d$ vertices $v_j \coloneqq \phi(j)$ as follows:
\begin{equation*}
P_i = \begin{cases}
(v_{i-2\kappa}, v_{i-2\kappa + 2}, \dots, v_i, \dots, v_{i + 2 \kappa}) & \text{if } d = 2\kappa + 1, \\
(v_{i - 2 \kappa + 1}, v_{i-2 \kappa + 3}, \dots, v_{i-1}, v_{i+1}, \dots, v_{i + 2 \kappa - 1}) & \text{if } d = 2 \kappa.
\end{cases}
\end{equation*}

\begin{definition}
	The \emph{short-diagonal pentagram map} $T$ takes each $v_i$ to the intersection of $d$ consecutive short-diagonal hyperplanes around $v_i$. Explicitly,
	\begin{equation*}
	T(v_i) = \begin{cases}
	P_{i - \kappa} \cap P_{i-\kappa + 1} \cap \dots \cap P_i \cap \dots \cap P_{i + \kappa} & \text{if } d = 2 \kappa + 1, \\
	P_{i - \kappa + 1} \cap P_{i-\kappa + 2} \cap \dots \cap P_i \cap \dots \cap P_{i + \kappa} & \text{if } d = 2 \kappa.
	\end{cases}
	\end{equation*}
	(When $d = 2$, this agrees with the construction of the pentagram map given earlier.)
\end{definition}

To describe a set of coordinates for $\P_n$, assume for technical reasons that $\gcd(n, d+1) = 1$. Given $\phi \in \P_n$, there exists a lift of the vertices $v_i = \phi(i) \in \RP^d$ to vectors $V_i \in \R^{d+1}$ such that
\begin{equation*}
	\det (V_i, V_{i+1}, \dots, V_{i+d}) = 1
\end{equation*}
for all $i \in \Z$; these vectors satisfy difference equations
\begin{equation} \label{eq:differenceeq}
	V_{i+d+1} = a_{i,d} V_{i+d} + a_{i, d-1} V_{i+d-1} + \dots + a_{i,1} V_{i+1} + (-1)^d V_i, \quad i \in \Z,
\end{equation}
where the coefficients $a_{i,j}$ are $n$-periodic in the first index.

The continuous limit as $n \to \infty$ of a twisted $n$-gon with monodromy $M$ has always been thought of as a smooth curve $\gamma : \R \to \RP^{d}$  with monodromy (i.e., such that $\gamma(x+ 2\pi) = M \cdot \gamma(x)$ for all~$x$). This interpretation will be justified in the next section. 
The  assumption that the points of the $n$-gon lie in general position translates to the condition that $\gamma$ is nondegenerate, i.e., the vectors $\gamma'(x), \dots, \gamma^{(d)}(x)$ are linearly independent for every $x \in \R$ in some (hence any) affine chart. 
Such a curve $\gamma$ has a lift $\Gamma$ to $\R^{d+1}$ such that 
\begin{equation*}\label{eq:lift_condition}
	\det (\Gamma, \Gamma', \dots, \Gamma^{(d)})(x) = 1
\end{equation*}
for all $x \in \R$; this lift then satisfies the differential equation
\begin{equation} \label{eq:differentialeq}
	\Gamma^{(d+1)} (x) + u_{d-1}(x) \Gamma^{(d-1)}(x) + \dots + u_1(x) \Gamma'(x) + u_0(x) \Gamma(x) = 0,
\end{equation}
for some $2\pi$-periodic functions $u_i$.
Thus to each such curve $\gamma$ we can associate a corresponding differential operator $L \coloneqq \partial^{d+1} + u_{d-1} \partial^{d-1} +\dots +  u_0$ (where $\partial = \frac{d}{dx}$ and $u_i \in C^\infty(S^1)$).

Given a curve $\gamma$ as above, define a family of hyperplanes $P_\epsilon(x)$ passing through $d$ points
\begin{equation*}
	P_\epsilon(x) = \begin{cases}
	(\gamma(x- \kappa \epsilon), \dots, \gamma(x), \dots, \gamma(x + \kappa \epsilon)) & \text{if } d = 2 \kappa +1, \\
	(\gamma(x- (2 \kappa - 1)\epsilon), \gamma(x - (2 \kappa -3 )\epsilon), \dots, \gamma(x + (2\kappa - 1)\epsilon)) & \text{if } d = 2 \kappa.
	\end{cases}
\end{equation*}
For fixed $\epsilon$, let $\gamma_\epsilon$ be the envelope curve of the hyperplanes $P_\epsilon(x)$, meaning that $\gamma_\epsilon(x)$ lies on $P_\epsilon(x)$ and the vectors $\gamma_\epsilon'(x), \dots, \gamma_\epsilon^{(d-1)}(x)$ span $P_\epsilon(x)$ for each $x$. 

\begin{definition}
    The envelope construction of the continuous limit of the pentagram map is given by the evolution of $\gamma$ in the direction of $\gamma_\epsilon$. 
    More precisely, let
    \begin{equation*}
    	L_\epsilon \coloneqq \partial^{d+1} + u_{d-1, \epsilon} \partial^{d-1} + \dots  + u_{0,\epsilon} 
    \end{equation*}
    be the differential operator corresponding to $\gamma_\epsilon$. 
    The coefficients $u_{i, \epsilon}\in C^\infty(S^1)$ have $\epsilon$-expansions of the form
    \begin{equation*}
    	u_{i,\epsilon}(x) = u_i(x) + \epsilon^2 w_i(x) + \O(\epsilon^4), \quad i = 0, \dots, d-1,
    \end{equation*}
    due to the symmetry $\epsilon \to -\epsilon$.
    Regarding $\epsilon^2$ as time, we define the continuous limit of the pentagram map by the evolution equations $du_i/dt = w_i$ for $i = 0, \dots, d-1$.
\end{definition}

As proved in Theorems 4.3 and 4.5 of~\cite{Khesin:2013}, the system $du_i/dt = w_i$ in any dimension $d$ corresponds to the $(2,d+1)$-KdV equation (defined below). For $d=2$, this is the classical Boussinesq equation (see Theorem 5 of~\cite{Ovsienko:2010}).

\subsection{The KdV hierarchy}

We give a brief description of the KdV hierarchy, primarily following~\cite{Adler:1978} (but in less generality).

Define the algebra $\mathcal{A}$ of formal pseudodifferential operators with periodic coefficients, consisting of formal series
\begin{equation*}
	B = \sum_{i=-\infty}^N b_i \partial^i, \quad b_i \in C^\infty(S^1),
\end{equation*}
with multiplication defined by
\begin{equation*}
	\partial^k \circ b = \sum_{n \geq 0} \binom{k}{n} b^{(n)} \partial^{k-n}, \quad k \in \Z,\, b \in C^\infty(S^1).
\end{equation*}
(Note that the  binomial coefficient $\binom{k}{n} \coloneqq \frac{k(k-1)\dotsm (k-n + 1)}{n!}$ also makes sense for $k < 0$.) Let $\DO \subseteq \mathcal{A}$ denote the subalgebra of differential operators.

Fix a differential operator of the form
\begin{equation*}
L = \partial^{d+1} + u_{d-1} \partial^{d-1} +\dots +  u_0,
\end{equation*}
where the $u_i$ are periodic functions. Define the root $L^{1/d+1}$ as the unique pseudodifferential operator of the form $\partial + \sum_{i < 0} b_i \partial^i$ satisfying $(L^{1/d+1})^{d+1}=L$ (the coefficients $b_i$ may be computed from this property). We may then define rational roots as $L^{m/d+1} \coloneqq (L^{1/d+1})^m$ for $m\in \Z_{\geq 0}$.

For any such $m$, let $Q_m \coloneqq (L^{m/d+1})_+$ denote the purely differential part of $L^{m/d+1}$ (i.e., its projection to $\DO$). 
In particular, one has $Q_2 = \partial^2 + \frac{2}{d+1} u_{d-1}$. 
Define an evolution equation on $L$ (that is, on its coefficients $u_i$) as follows:

\begin{definition}
	The \emph{$(m, d+1)$-KdV equation} is the evolution equation
	\begin{equation*}
		\frac{d}{dt}L = [Q_m, L].
	\end{equation*}
\end{definition}
\begin{remark}\label{rem:kdv-lax}
	The KdV equations can be described as Hamiltonian equations on $\DO$ (see Theorem~4 of~\cite{Adler:1978}); these systems are known to be completely integrable. 
	Note also that the above equation is written in Lax form; we discuss this further in Section~\ref{sec:lax_forms}.
\end{remark}

We will frequently use the following fact:

\begin{prop}[Proposition~10.1 of \cite{Khesin:2013}] \label{prop:gamma_evolution_to_kdv}
    Suppose $\Gamma$ is a nondegenerate curve in $\R^{d+1}$ satisfying $L \Gamma = 0$. Then the evolution $d\Gamma/dt = Q_m \Gamma$ implies the $(m, d+1)$-KdV equation $dL/dt =[Q_m, L]$.
\end{prop}


\section{Continuous limits of higher-dimensional pentagram maps} \label{sec:ctslimitcalc}

In this section, we provide a rigorous and direct treatment of the continuous limit for the short-diagonal pentagram map and several generalizations.
The standard construction of the continuous limit involves two main components: passing from twisted $n$-gons to nondegenerate curves with monodromy, and using the envelope curve as an analogue of the short-diagonal map. 
We will justify both of these components by employing quantum calculus (see~\cite{Kac:2002}) and explicitly comparing them to limits of their discrete counterparts.
Furthermore, we introduce the $\chi$-pentagram map, which generalizes the construction to a broad class of maps involving intersections of different subspaces.

\subsection{Kinematics: Limit of the phase spaces} 

\label{sub:limits_of_difference_equations}

Here, we show that the differential equation \eqref{eq:differentialeq} corresponding to a nondegenerate curve can be recovered by taking a limit of difference equations corresponding to twisted $n$-gons. We begin by fixing a nondegenerate curve $\gamma: \R \to \RP^d$ with monodromy $M$. Let $\Gamma$ be its lift to $\R^{d+1}$ such that $\det (\Gamma, \Gamma', \dots, \Gamma^{(d)}) \equiv 1$. This lift satisfies \eqref{eq:differentialeq}; that is, we have
\begin{equation*} 
	\Gamma^{(d+1)}+ u_{d-1}\Gamma^{(d-1)} +  \dots + u_0 \Gamma\equiv 0
\end{equation*}
for some $2\pi$-periodic functions $u_i$.

Following~\cite{Ovsienko:2010}, we discretize the curve by fixing $x \in \R$ and some small $\epsilon > 0$, and setting $v_i \coloneqq \gamma(x + i \epsilon)$.
(In particular, for $\epsilon = 2 \pi /n$, the $(v_i)$ are the vertices of a twisted $n$-gon ``converging'' to $\gamma$ as $n \to \infty$.)

Lift each $v_i$ to the vector $\tilde{V}_i \coloneqq \Gamma(x + i \epsilon)$ in $\R^{d+1}$. For small $\epsilon$, nondegeneracy of $\Gamma$ implies that $\tilde{V}_0, \dots, \tilde{V}_d$ are linearly independent for all $x$.    
We can therefore write
\begin{equation*}
	\tilde{V}_{d+1} = \tilde{a}_d \tilde{V}_d + \dots +  \tilde{a}_0 \tilde{V}_0
\end{equation*}
for some coefficients $\tilde{a}_i = \tilde{a}_i(x,\epsilon)$, periodic in $x$. More explicitly, 
\begin{equation} \label{eq:diffeq-v1}
	\Gamma (x + (d+1) \epsilon) = \tilde{a}_d(x, \epsilon) \Gamma(x + d \epsilon) + \dots + \tilde{a}_0(x, \epsilon) \Gamma(x).
\end{equation}

\begin{remark} \label{rem:vertices_vs_curve_lift}
	We emphasize that the $\tilde{V}_i$ are obtained via the canonical lift of the curve $\gamma \to \Gamma$, as opposed to the canonical lift of the vertices $v_i$. 
	This differs from the setting of~\cite{Khesin:2013}, since the $\tilde{V}_i$ in general do not satisfy the difference equation \eqref{eq:differenceeq}; that is, we do not have $\tilde{a}_0 \equiv (-1)^d$. However, we will later see that $\tilde{a}_0(x,\epsilon) \to (-1)^d$ for each $x \in \R$ as $\epsilon \to 0$.
\end{remark}

It will be convenient for what follows to define a difference operator $\Delta_\epsilon$ by
\begin{equation*}
	\Delta_\epsilon \Gamma(x) \coloneqq \Gamma(x+\epsilon) - \Gamma(x).
\end{equation*}
(cf.~the $h$-derivative defined in~\cite{Kac:2002}, where one can find many other properties of this operator.)
Rewrite \eqref{eq:diffeq-v1} as
\begin{equation}\label{eq:diffeq-v2}
	(\Delta_\epsilon^{d+1} + A_{d} \Delta_\epsilon^{d} + A_{d-1} \Delta_\epsilon^{d-1} + \dots + A_0) \Gamma = 0
\end{equation}
for certain coefficients $A_i = A_i(x, \epsilon)$,  again periodic in $x$.

Note that one can rearrange the difference equation \eqref{eq:diffeq-v2} back to its original form \eqref{eq:diffeq-v1} using the expansion
\begin{equation*} \label{eq:difference-operator-expansion}
    \Delta_\epsilon^k \Gamma(x) = \sum_{i=0}^k (-1)^{k-i} \binom{k}{i} \Gamma(x + i \epsilon).
\end{equation*}
By equating the corresponding coefficients, we find that
\begin{equation} \label{eq:a_i vs A_i}
    \tilde{a}_i = \sum_{k=i}^{d+1} (-1)^{k-i+1} \binom{k}{i} A_k \quad (\text{where } A_{d+1} \coloneqq 1).
\end{equation}

\begin{prop}\label{prop:limit-of-difference-eq}
    The differential equation \eqref{eq:differentialeq} is the quantum calculus limit of the difference equation \eqref{eq:diffeq-v2}. More precisely, 
    \begin{equation*} 
        \Delta_\epsilon^{d+1} \Gamma + A_{d} \Delta_\epsilon^{d}\Gamma + \dots + A_0\Gamma
        = \epsilon^{d+1} (\Gamma^{(d+1)} + u_{d-1}\Gamma^{(d-1)} +  \dots + u_0 \Gamma) + \O(\epsilon^{d+2}),
    \end{equation*}
i.e.~Equation~\eqref{eq:differentialeq} is the lowest degree term of Equation~\eqref{eq:diffeq-v2} as $\epsilon \to 0$. 
\end{prop}

Before we calculate the coefficients $A_i$ in order to prove the proposition, we mention some elementary technical properties of $\Delta_\epsilon$.

\begin{lemma}\label{lem:Taylor}
	Expanded as a series in $\epsilon$, we have
	\begin{equation*}
		\Delta_\epsilon^k \Gamma(x) = \epsilon^k \,\Gamma^{(k)}(x)  + \O(\epsilon^{k+1})
	\end{equation*}
	for any $x \in \R$ and $k \in \Z_{\geq 0}$. 
    The higher-order terms of the expansion are of the form $\epsilon^\ell c_{\ell, k} \Gamma^{(\ell)}(x)$ for some constants $c_{\ell, k}$.
\end{lemma}
\begin{proof}
The $k=1$ case is clear by Taylor expansion, and the general case follows by induction.
\end{proof}

In other words,
\begin{equation} \label{eq:differencelimit}
    \lim_{\epsilon \to 0} \frac{\Delta_\epsilon^k \Gamma(x)}{\epsilon^k}  = 
\Gamma^{(k)}(x).
\end{equation}

Next we calculate the coefficients $A_{i}$ of \eqref{eq:diffeq-v2}.

\begin{lemma}\label{lem:A_i expansions}
	Expanding each coefficient $A_i = A_i(x,\epsilon)$ as a series in $\epsilon$, we have
	\begin{equation*}
		A_i = \epsilon^{d + 1 -i}  \, u_i  + \O(\epsilon^{d+2-i}) 
	\end{equation*}
	for $i=0, \dots, d$ (with the convention $u_d \coloneqq 0$). 
\end{lemma}
\begin{proof}
	This follows from the Taylor expansion of \eqref{eq:diffeq-v2} up to the $\epsilon^{d+1}$ order terms, which is straightforward with the use of Lemma~\ref{lem:Taylor}. 
\end{proof}

\begin{remark}\label{rmk:A_i refinement}
    By examining the $\epsilon^{d+2}$ term in the Taylor expansion, one can see that 
    \begin{equation*}
        A_d = \epsilon^2 u_{d-1}  + \O(\epsilon^3).
    \end{equation*}
    This refinement is not required for the proof of Proposition~\ref{prop:limit-of-difference-eq}, but will be useful later.
    More generally, the higher-order terms of the $A_i$ are differential polynomials in the $u_j$.
\end{remark}

\begin{proof}[Proof of Proposition~\ref{prop:limit-of-difference-eq}]
   Expanding the coefficients $A_i$ in \eqref{eq:diffeq-v2} and dividing by $\epsilon^{d+1}$, we get
\begin{equation*}
	\frac{\Delta_\epsilon^{d+1} \Gamma(x)}{\epsilon^{d+1}} + (0 + \O(\epsilon)) \frac{\Delta_\epsilon^d \Gamma(x)}{\epsilon^{d}} + \left( u_{d-1} + \O(\epsilon) \right) \frac{\Delta_\epsilon^{d-1} \Gamma(x)}{\epsilon^{d-1}} + \dots + (u_0 + \O(\epsilon)) \Gamma(x) = 0. 
\end{equation*}
Letting $\epsilon \to 0$ and using \eqref{eq:differencelimit}, we recover the desired differential equation
\begin{equation*}
	\Gamma^{(d+1)}(x) + u_{d-1} \Gamma^{(d-1)}(x) + \dots + u_0 \Gamma(x) = 0. \qedhere
\end{equation*}
\end{proof}

\begin{remark}
To conclude this discussion, we give an expansion of the coefficients $\tilde{a}_i$.
By substituting the expansions of the $A_i$ from Lemma~\ref{lem:A_i expansions} and Remark~\ref{rmk:A_i refinement} into \eqref{eq:a_i vs A_i}, we obtain
	\begin{align*}
	    \tilde{a}_i(x,\epsilon) 
	    &= (-1)^{d-i}  \binom{d+1}{i} +  \epsilon^2 (-1)^{d-i+1} \binom{d-1}{i-1} u_{d-1}(x)   + \O(\epsilon^3).
	\end{align*}
In particular, one finds that $\tilde{a}_0(x,\epsilon) = (-1)^d + \O(\epsilon^3)$, as we alluded to in Remark~\ref{rem:vertices_vs_curve_lift}.

    In the case $d =2$, the coefficients that we obtain agree with those given in Lemma 6.9 of~\cite{Ovsienko:2010}, up to the second-order terms. 
    However, the third-order terms differ because of the assumption that $\tilde{a}_0 \equiv 1$ in~\cite{Ovsienko:2010}. 
\end{remark}

\subsection{Dynamics: Limit of the short-diagonal pentagram map} 
\label{sub:limit_of_short_diag_pentagram_map}

Next, we wish to justify the use of the envelope curve in the continuous limit. To this end, we show that applying the pentagram map directly to points on a curve yields the same evolution (after a suitable reparametrization) as the envelope construction, i.e., it also gives the $(2, d+1)$-KdV equation. 

We define a continuous analogue $T_\epsilon$ of the short-diagonal map on curves as follows.
Using the same notation as before, we begin by fixing a nondegenerate curve $\gamma$ in $\RP^d$, which we discretize for fixed $x$ and small $\epsilon$  by setting $v_i \coloneqq \gamma(x + i \epsilon)$. 
Let $T_\epsilon(\gamma)(x)$ denote the image of the point $\gamma(x)$ under the pentagram map, using the vertices $(v_i)$.
We can thus obtain a new curve $\tilde{\gamma}_\epsilon$ by

\begin{equation}\label{eq:shifted_short_diag}
    \tilde{\gamma}_\epsilon(x) = \begin{cases}
        T_\epsilon(\gamma)(x) & \text{if $d$ is odd}, \\
         T_\epsilon(\gamma)(x- \frac{\epsilon}{2}) & \text{if $d$ is even.} 
    \end{cases}
\end{equation}
Let $\tilde{\Gamma}_\epsilon$ denote the lift of $\tilde{\gamma}_\epsilon$ to $\R^{d+1}$ with normalization $\det (\tilde{\Gamma}_\epsilon, \tilde{\Gamma}_\epsilon',\dots, \tilde{\Gamma}_\epsilon^{(d)}) \equiv 1$.

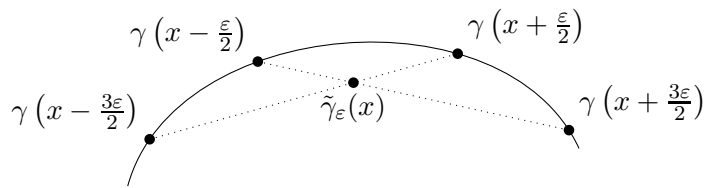
\begin{figure}[ht] 
	\begin{centering}
	\begin{tikzpicture}
    	\draw (0,0) to [bend left=70] 
        coordinate[pos=0.1] (v_{i-1})
    	coordinate[pos=0.35] (v_i) 
    	coordinate[pos=0.7] (v_{i+1})
    	coordinate[pos=0.96] (v_{i+2}) (6,0.5);
    	\fill[black] (v_i) circle (2pt);
    	\fill[black] (v_{i-1}) circle (2pt);
    	\fill[black] (v_{i+1}) circle (2pt);
    	\fill[black] (v_{i+2}) circle (2pt);
    	\draw[dotted] (v_i) -- (v_{i+2});
    	\draw[dotted] (v_{i-1}) -- (v_{i+1});
    	\coordinate (u) at (intersection of {v_{i-1}--v_{i+1}} and {v_i--v_{i+2}});
    	\fill[black] (u) circle (2pt);
        \node at (v_i) [above left] {$\gamma\left(x-\frac{\epsilon}{2}\right)$};
        \node at (v_{i+1}) [above right] {$\gamma\left(x+\frac{\epsilon}{2}\right)$};
        \node at (v_{i+2}) [above right] {$\gamma\left(x+\frac{3\epsilon}{2}\right)$};
        \node at (v_{i-1}) [above left] {$\gamma \left(x - \frac{3 \epsilon}{2} \right)$};
        \node at (u) [below] {$\tilde{\gamma}_\epsilon(x)$};
    \end{tikzpicture}
\end{centering}
	\caption{The short-diagonal map on a discretized curve ($d=2$)}
	\label{fig:discretized curve}
\end{figure}

\begin{remark}\label{rem:short_diagonal_shift}
	Figure~\ref{fig:discretized curve} gives motivation for the shift by $\epsilon/2$ when $d$ is even: the point $T_\epsilon(\gamma)(x- \frac{\epsilon}{2}) = \tilde{\gamma}_\epsilon(x)$ is centred near $\gamma(x)$, and there is symmetry $\epsilon \to -\epsilon$.
    This shift turns out to be necessary to ensure that the expansion of $\tilde{\Gamma}_\epsilon$ as a series in $\epsilon$ has no linear term.
\end{remark}

Our main result in this section is the following analogue of Theorem 4.3 in~\cite{Khesin:2013}.

\begin{theorem}\label{thm:directctslimit}
     The curve $\tilde{\Gamma}_\epsilon$ has the expansion $\tilde{\Gamma}_\epsilon = \Gamma + \epsilon^2 \tilde{C}_d \cdot Q_2 \Gamma + \O(\epsilon^4)$, i.e.
    \begin{equation*}
        \tilde{\Gamma}_\epsilon(x) = \Gamma(x) + \epsilon^2 \tilde{C}_d \left( \Gamma''(x) + \frac{2}{d+1} u_{d-1}(x) \Gamma(x) \right) + \O(\epsilon^4)
    \end{equation*}
    as $\epsilon \to 0$, for some nonzero constant $\tilde{C}_d$.
\end{theorem}
\begin{remark}
    The constant $\tilde{C}_d$ is not the same as the constant $C_d$ from~\cite{Khesin:2013} (which uses the envelope construction of the continuous limit). 
    However, this only affects the time parametrization; both cases lead to the evolution $d \Gamma/dt = Q_2 \Gamma$ (where $Q_2\coloneqq (L^{2/d+1})_+ = \partial^2 + \frac{2}{d+1} u_{d-1}$), corresponding to the $(2,d+1)$-KdV equation $dL/dt = [Q_2, L]$.
    Consequently this justifies the use of the envelope curve $\gamma_\epsilon$ in the definition of the continuous limit, 
    i.e., shows that the envelope construction is equivalent to that involving direct application of the pentagram map and taking the quantum calculus limit.
\end{remark}

Before proving the theorem, we give a more explicit description of $\tilde{\Gamma}_\epsilon$. 
For simplicity in notation, we will assume that $d = 2 \kappa + 1$ is odd. (The same argument will work for even~$d$.)
By the definition of the short-diagonal pentagram map given in Section~\ref{sub:short_diag_background}, the point $\tilde{\gamma}_\epsilon(x)= T_\epsilon(\gamma)(x)$ lies on the planes 
\begin{equation*}
    P_j \coloneqq (\gamma(x + (j-2\kappa)\epsilon), \gamma(x + (j-2\kappa+2)\epsilon), \dots, \gamma(x + (j+2\kappa)\epsilon))
\end{equation*}
for $j = - \kappa, -\kappa + 1, \dots, \kappa$. Lifting to $\R^{d+1}$, we obtain the corresponding conditions
\begin{equation}\label{eq:odd-condition1}
    0 = \det(\tilde{\Gamma}_\epsilon(x), \Gamma(x + (j-2\kappa)\epsilon), \Gamma(x + (j - 2 \kappa + 2)\epsilon), \dots, \Gamma(x + (j + 2 \kappa)\epsilon)).
\end{equation}

Observe that if we replace $\epsilon$ by $-\epsilon$, then $\tilde{\Gamma}_{-\epsilon}(x)$ satisfies the same defining equations as $\tilde{\Gamma}_\epsilon(x)$. (Indeed, the condition \eqref{eq:odd-condition1} for each value of $j$ switches with that for $-j$.)
It follows that the $\epsilon$-expansion of $\tilde{\Gamma}_\epsilon$ has only even powers of $\epsilon$; we will write 
\begin{equation*}
    \tilde{\Gamma}_\epsilon = \Gamma + \epsilon^2 B + \O(\epsilon^4)
\end{equation*} 
for some function $B = B(x)$.

In order to calculate $B$, it will be useful to rewrite \eqref{eq:odd-condition1} as follows.
For each equation indexed by~$j=-\kappa, -\kappa+1,  \ldots, \kappa$ in~\eqref{eq:odd-condition1}, replace $x$ by $x + (\kappa - j)\epsilon$ in order to obtain the system
\begin{equation*}
    0 = \det(\tilde{\Gamma}_\epsilon(x + i\epsilon), \Gamma(x -\kappa \epsilon), \Gamma(x - (\kappa - 2)\epsilon), \dots, \Gamma(x + 3 \kappa\epsilon))
\end{equation*}
for $i = 0,1, \dots, d-1$. 
By linearity in the first column it follows that
\begin{equation}\label{eq:odd-condition2}
    0 = \det(\Delta_\epsilon^i \tilde{\Gamma}_\epsilon(x), \Gamma(x -\kappa \epsilon), \Gamma(x - (\kappa - 2)\epsilon), \dots, \Gamma(x + 3 \kappa\epsilon))
\end{equation}
for each such $i$, where $\Delta_\epsilon$ is the difference operator defined in Section~\ref{sub:limits_of_difference_equations}.

We now complete the proof of the main result (cf.~\cite{Khesin:2013}).

\begin{proof}[Proof of Theorem~\ref{thm:directctslimit}]
    Start by writing $B = b_0 \Gamma + \dots + b_{d} \Gamma^{(d)}$ for some coefficients $b_i = b_i(x)$. (Note that the vectors $\Gamma(x), \dots, \Gamma^{(d)}(x)$ form a basis for each $x$, since $\det(\Gamma, \dots, \Gamma^{(d)}) \equiv 1$.)

    We will expand \eqref{eq:odd-condition2} in $\epsilon$ for $0 \leq i \leq d-1$ and examine the lowest  terms. By Lemma~\ref{lem:Taylor}, we know that
    \begin{equation*}
        \Delta_\epsilon^i \tilde{\Gamma}_\epsilon = \epsilon^i \tilde{\Gamma}_\epsilon^{(i)} + \epsilon^{i+1} c_i \tilde{\Gamma}_\epsilon^{(i+1)} + \epsilon^{i+2} d_i \tilde{\Gamma}_\epsilon^{(i+2)} + \dots
    \end{equation*}
    for some constants $c_i, d_i$. Since $\tilde{\Gamma}_\epsilon = \Gamma + \epsilon^2 B + \O(\epsilon^4)$, it follows that
    \begin{equation}\label{eq:tilde-diff-expansion}
        \Delta_\epsilon^i \tilde{\Gamma}_\epsilon = \epsilon^i \Gamma^{(i)} + \epsilon^{i+1} c_i \Gamma^{(i+1)} + \epsilon^{i+2} (B^{(i)} + d_i \Gamma^{(i+2)}) + \O(e^{i+3}).
    \end{equation}

    For $0 \leq i < d-2$, note that the $\Gamma^{(i)}, \Gamma^{(i+1)}$, and $\Gamma^{(i+2)}$ terms in \eqref{eq:tilde-diff-expansion} are all killed by $\det(\cdot, \Gamma, \Gamma', \dots, \Gamma^{(d-1)})$. 
    Hence the $\epsilon^{(i+2) + 0 + 1 + \dots + (d-1)} = \epsilon^{(i+2) + \frac{d(d-1)}{2}}$ term of \eqref{eq:odd-condition2} yields
    \begin{equation*}
        \det(\Gamma, \dots, \Gamma^{(d-1)}, B^{(i)}) = 0,
    \end{equation*}
     which implies that $b_{d-i} \equiv 0$.

     For $i = d-2$, the $\epsilon^{d + \frac{d(d-1)}{2}}$ term of \eqref{eq:odd-condition2} implies that
     \begin{equation*}
         \det(\Gamma, \dots, \Gamma^{(d-1)}, B^{(d-2)}) = \tilde{C}_d \det(\Gamma, \dots, \Gamma^{(d-1)}, \Gamma^{(d)}) = \tilde{C}_d
     \end{equation*}
     for some nonzero constant $\tilde{C}_d$, and hence $b_2 \equiv \tilde{C}_d$.

     Finally, for $i = d-1$, the $\epsilon^{(d + 1) + \frac{d(d-1)}{2}}$ term of \eqref{eq:odd-condition2} again yields
     \begin{equation*}
         \det(\Gamma, \dots, \Gamma^{(d-1)}, B^{(d-2)}) = 0.
     \end{equation*}
     (Here we use the fact that $\det(\Gamma, \dots, \Gamma^{(d-1)}, \Gamma^{(d+1)}) \equiv 0$.) It follows that $b_{d-1} \equiv 0$.

     Thus we have shown $\tilde{\Gamma}_\epsilon = b_0 \Gamma + \tilde{C}_d \Gamma''$. The rest of the proof is exactly the same as the proof of Theorem 4.3 in~\cite{Khesin:2013}. By expanding the normalization condition $\det(\tilde{\Gamma}_\epsilon, \tilde{\Gamma}_\epsilon', \dots, \tilde{\Gamma}_\epsilon^{(d)}) \equiv 1$, one can show that $b_0 = \tilde{C}_d \cdot \frac{2}{d+1} u_{d-1}$, which gives the desired result.
\end{proof}

\begin{remark}
    The short-diagonal pentagram map can be generalized by choosing different vertices to make up each hyperplane, and by choosing different hyperplanes to intersect (see the \emph{generalized pentagram map} $T_{I,J}$ defined in~\cite{Khesin:2016}).
    As long as the hyperplanes are consecutive, the same discretization argument as above gives a direct construction of the the continuous limit.
    The choice of points making up each hyperplane only affects the constants present in the calculations, and in nondegenerate cases one again obtains the $(2,d+1)$-KdV equation.
    In Example~\ref{ex:dual_dented} below, we consider a case involving nonconsecutive hyperplanes, requiring different techniques.
\end{remark}

\subsection{Limit of the \texorpdfstring{$\chi$-pentagram}{χ-pentagram} map}\label{sub:limit_of_chi_pent}

The short-diagonal map and its generalizations discussed thus far have all been found to yield the $(2,d+1)$-KdV equation as their continuous limit. 
In this section, we seek to generalize this result to a broader class of geometric constructions and provide some justification for the appearance of this equation. 
To this end, we introduce the \textit{$\chi$-pentagram map} (describing a ``\textbf{C}ontinuous \textbf{H}igher \textbf{I}ntersection'' generalization of the pentagram map!) below.

\begin{definition}
    Fix a set $\chi=\{\{p_{1,0},\dots, p_{1,q_1}\},\dots,\{p_{r,0},\dots, p_{r,q_r}\} \}$ with $p_{i,j}\in \R$ such that $p_{i,j}\neq p_{i,k}$, and $q_i\in \Z_{\geq 1}$ such that $(d-q_1)+\dots+(d-q_r)=d$. For  $\epsilon > 0$ and $\gamma\subset \RP^d$ a nondegenerate curve, consider the subspaces
    \begin{equation*}
        P_\e^i(x)=\big(\gamma(x+p_{i,0}\e),\dots, \gamma(x+p_{i,q_i}\e)\big)
    \end{equation*}
    for $1\leq i\leq r$. The \emph{$\chi$-pentagram map} sends $\gamma$ to the curve $T_\e^\chi(\gamma)$, where
    \begin{equation*}
        T_\e^\chi(\gamma)(x)=\bigcap_{i=1}^r P_\e^i(x).
    \end{equation*}
\end{definition}

\begin{remark}
    This construction encompasses the definition of the short-diagonal map on a discretized curve from Section~\ref{sub:limit_of_short_diag_pentagram_map}, as well as the generalized pentagram map from~\cite{Khesin:2016}, hence the difficulty in finding a name. When the parameters $p_{i,j}$ are not restricted to integer values, the map $T_\e^\chi$ does not have a polygonal analogue (unlike the case in Section~\ref{sub:limit_of_short_diag_pentagram_map}).
\end{remark}

As before, we lift $\gamma$ to $\Gamma \subset \R^{d+1}$ with the normalization $\det(\Gamma, \Gamma',\dots, \Gamma^{(d)}) \equiv 1$, satisfying $L \Gamma = 0$ for $L = \partial^{d+1} + u_{d-1} \partial^{d-1} + \dots + u_0$.
Similarly, let $\ge\coloneqq\ge^\chi$ denote the normalized lift of $T_\e^\chi(\gamma)$, and expand $\ge$ in $\e$ as 
\begin{equation}\label{eq:ge_expansion}
    \ge=\Gamma+\e G_1\Gamma+\e^2 G_2\Gamma+\dots, \quad \tn{where } G_i=\sum_{j=0}^d\alpha_{i,j}\partial^j 
\end{equation}
for some coefficients $\alpha_{i,j} = \alpha_{i,j}(x)$. 
By convention, we set $G_0 = 1$ (i.e.~$\alpha_{0,j} = \delta_{0,j}$).
We now turn to the problem of computing the limit of the described evolution as $\e\to0$.

In order to generalize our earlier arguments involving determinants, we consider the wedge product on $\bigwedge \R^{d+1}$. 
Note that the wedges of the derivatives $\Gamma^{(i)}(x)$ form a basis of $\bigwedge \R^{d+1}$ for each $x \in \R$, by nondegeneracy of $\Gamma$.
From the definition of the $\chi$-pentagram map, the coplanarity conditions are
\begin{equation} \label{eq:chi_wedge_condition}
    \ge(x) \wedge \Gamma(x+p_{i,0}\e)\wedge\dots\wedge \Gamma(x+p_{i,q_i}\e)=0
\end{equation}
for $1\leq i\leq r$ (corresponding to each space $P_\e^i$).

\begin{example}[Short-diagonal map]
    For odd $d = 2 \kappa + 1$, say, and
    \begin{equation*}
        \chi = \{ \{i-2 \kappa, i-2 \kappa + 2 , \dots, i + 2 \kappa\} : - \kappa \leq i \leq \kappa \},
    \end{equation*}
    we recover the construction \eqref{eq:shifted_short_diag} of the short-diagonal pentagram map on a discretized curve, and the coplanarity conditions \eqref{eq:chi_wedge_condition} are equivalent to the earlier determinant conditions \eqref{eq:odd-condition1}.
\end{example}

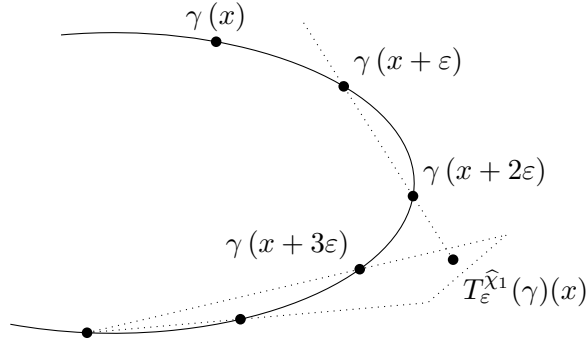
\begin{figure}[ht]
\begin{centering}
\begin{tikzpicture}

    \draw (1,0) arc (100:-110:4cm and 2cm) 
        coordinate[pos=1/7] (v_0)
        coordinate[pos=2/7] (v_1)
        coordinate[pos=7/14] (v_2)
        coordinate[pos=9/14] (v_3)
        coordinate[pos=11/14] (v_4)
        coordinate[pos=13/14] (v_5);

    \fill[black] (v_0) circle (2pt);
    \fill[black] (v_1) circle (2pt);
    \fill[black] (v_2) circle (2pt);
    \fill[black] (v_3) circle (2pt);
    \fill[black] (v_4) circle (2pt);
    \fill[black] (v_5) circle (2pt);




    \draw[dotted, shorten >=-1cm, shorten <=-1cm] (v_1) -- (v_2);
    \draw[dotted, shorten >=-2cm] (v_5) -- (v_3);
    \draw[dotted, shorten >=-2.5cm] (v_5) -- (v_4);

    \coordinate (l) at ($(v_2)!-1cm!(v_1)$);
    \coordinate (p_1) at ($(v_3)!-2cm!(v_5)$);
    \coordinate (p_2) at ($(v_4)!-2.5cm!(v_5)$);

    \draw[dotted] (p_1) -- (p_2); 
   
    \fill[black] (l) circle (2pt);

    \node at (v_0) [above] {$\gamma\left(x\right)$};
    \node at (v_1) [above right] {$\gamma\left(x + \epsilon \right)$};
    \node at (v_2) [above right] {$\gamma\left(x + 2 \epsilon\right)$};
    \node at (v_3) [above left] {$\gamma\left(x + 3 \epsilon\right)$};
    
    \node at (l) [below right] {$T_\epsilon^{\widehat{\chi}_1}(\gamma)(x)$};
\end{tikzpicture}
\end{centering}
\caption{The dual dented map on a discretized curve ($d=3, s=1$)}
\label{fig:dual_dented}
\end{figure}

\begin{example}[Dual dented map]\label{ex:dual_dented}
    Fix $s \in \{1, \dots, d-1\}$. Given a generic twisted $n$-gon in $\RP^d$ with vertices $(v_j)$, consider for each $i \in \Z$ the hyperplane $P_{i} = (v_i, v_{i+1}, \dots, v_{i+d-1})$ through $d$~consecutive vertices. The \emph{dual dented pentagram map} $\widehat{T}_s$ (defined in~\cite{Khesin:2016}) takes each $v_i$ to the intersection of $d$ planes
    \begin{equation*}
        \widehat{T}_s (v_i) = P_i \cap P_{i+1} \cap \dots \cap P_{i+d-s-1} \cap P_{i+d-s+1} \cap  \dots \cap P_{i+d}
    \end{equation*}
    (i.e., skipping the plane $P_{i+d-s}$).
    The dynamics of this map is studied in~\cite{Khesin:2016}. 
    Following the approach of Section~\ref{sub:limit_of_short_diag_pentagram_map}, the analogue of $\widehat{T}_s$ on curves is given by the map $T_\e^{\widehat{\chi}_s}$ with
    \begin{equation*}
        \widehat{\chi}_s \coloneqq \{ \{i, i+1, \dots, i + d-1  \} : 0 \leq i \leq d, \,  i \neq d-s \}.
    \end{equation*}
    (According to~\cite{Khesin:2016}, we may equivalently take
    \begin{equation*} 
        \widehat{\chi}_s^{\,\prime} \coloneqq\l\{\{d-s-1,\dots, d-1\},\{d,\dots, 2d-s\}\r\},
    \end{equation*}
    describing the dual dented map as an intersection of two subspaces of complementary dimensions $s$ and $d-s$; see Figure~\ref{fig:dual_dented}.) Note that the difference operator argument used to compute the continuous limit of the short-diagonal map in Section~\ref{sub:limit_of_short_diag_pentagram_map} is not applicable to $T_\epsilon^{\widehat{\chi}_s}$, as it required that all of the hyperplanes to be intersected were consecutive. 
    The continuous limit of the dual dented map is described in Proposition~\ref{prop:dual_dented_limit}.
\end{example}

\begin{definition}
    We say that a configuration $\chi$ is \emph{centralized} if $G_1=0$, i.e.~the $\e$-expansion of $\ge$ is of the form
    \begin{equation*}
        \ge = \Gamma + \epsilon^2 G_2 \Gamma + \O(\epsilon^3).
    \end{equation*}

\end{definition}

\begin{example}
    The short-diagonal map described in Section~\ref{sub:limit_of_short_diag_pentagram_map} is centralized (see Remark~\ref{rem:short_diagonal_shift}).
\end{example}

Centralization is an algebraic condition on the parameters $p_{i,j}$ of $\chi$. Although this condition is hard to explicitly describe in general, we will characterize certain important examples in Section~\ref{sub:centralizing_chi_map}.

\begin{theorem}\label{thm:chi_pent_cont_limit}
    For any configuration $\chi$, the $\e$-expansion $\ge = \Gamma + \epsilon G_1 \Gamma + \epsilon^2 G_2 \Gamma + \O(\epsilon^3)$ is given by
    \begin{gather}
        G_1=\alpha_{1,1}\partial, \label{eq:G_1_expansion}\\
        G_2=\alpha_{2,2}\l(\partial^2+\frac{2}{d+1} u_{d-1}  \r) -  \frac{\alpha_{1,1}^2}{d+1} u_{d-1}, \label{eq:G_2_expansion}
    \end{gather}
    where the coefficients $\alpha_{1,1}$ and $\alpha_{2,2}$ are constant in $x$. 
    
    In particular, if $\chi$ is centralized and $G_2 \neq 0$, then the continuous limit of the $\chi$-pentagram map defined by 
    \begin{equation*}
        \frac{d}{dt} \Gamma = G_2 \Gamma
    \end{equation*}
    corresponds to the $(2,d+1)$-KdV equation.
\end{theorem}

\begin{remark}
    If $G_1\neq 0$, the corresponding continuous limit defined by $d\Gamma/dt = G_1 \Gamma = \alpha_{1,1} \Gamma'$ leads to the $(1,d+1)$-KdV equation, see e.g.~\cite{Khesin:2013}.
    In Section~\ref{sec:kdv_realization}, we consider situations in which higher $G_i$ also vanish.
\end{remark}

While the proof of Theorem \ref{thm:chi_pent_cont_limit} presented below resembles that of Theorem 4.3 in~\cite{Khesin:2013} (and its generalizations in e.g.~\cite{Khesin:2016, Izosimov:2019}), we emphasize the broadness of our setting. We do not, for instance, require that $\chi$ consists of arithmetic progressions (as in~\cite{Izosimov:2019}).
Theorem~\ref{thm:chi_pent_cont_limit} therefore further reinforces the idea that the $(2,d+1)$-KdV equation is a robust limit for pentagram-type maps (cf.~Remark 4.2 of~\cite{Khesin:2013}).

We begin by analyzing the coefficients $\alpha_{i,j}$ from the expansion \eqref{eq:ge_expansion}.

\begin{lemma}\label{lem:higher_coef_zero}
	For each $0 \leq i \leq d$, the differential operator $G_i$ from the $\epsilon$-expansion of $\ge$ is of order $\leq i$; that is, the coefficients $\alpha_{i,j}=0$ for $j>i$.
\end{lemma}
\begin{proof}
Recall that nondegeneracy of $\Gamma$ implies that we can write
\begin{equation*}
     \ge(x)=c_0 \Gamma(x)+ c_1 \Gamma(x+\e) + \dots+c_d \Gamma(x+d\e)
 \end{equation*}
with coefficients $c_j=c_j(x,\epsilon)$. 
Expanding the right-hand side in $\epsilon$, we see that for $0 \leq i \leq d$ the $\epsilon^i$ term (namely $G_i \Gamma$) involves only the derivatives $\Gamma, \dots, \Gamma^{(i)}$.
\end{proof}

\begin{prop} \label{prop:chi_pent_coef_expansion}
    Each $\alpha_{i,j}$ is a differential polynomial in the functions $u_\ell$, with coefficients depending on $\chi$.
\end{prop}

\noindent{\it Proof.}
This statement is proved by the following induction.
Let $k \in \Z_{\geq 0}$, and suppose that we have obtained the coefficients $\alpha_{j,0}, \alpha_{j+1, 1}, \dots, \alpha_{j+d,d}$ for all $0 \leq j \leq k-1$, as differential polynomials in the functions $u_\ell$. We may then solve for $\alpha_{k,0}, \alpha_{k+1,1}, \dots, \alpha_{k+d,d}$ as follows.
    
Start with  the $\epsilon$-expansion of the normalization condition $\det(\ge, \ge', \dots, \ge^{(d)}) \equiv 1$, from which
one can solve for $\alpha_{k,0}$ as a differential polynomial in the previously-obtained coefficients.
Indeed, the $\epsilon^k$ term of the expansion yields an expression for $\alpha_{k,0}$ as a differential polynomial in the other coefficients of $G_0, G_1, \dots, G_{k}$, all of which are known.

Now the proposition will be proved along with the following technical lemma.

\begin{lemma}\label{lem:coef_computation}
    For each $k \in \Z_{\geq 0}$ the $\epsilon$-expansions of the coplanarity conditions \eqref{eq:chi_wedge_condition}
    yield  a system
        \begin{equation} \label{eq:M_system} 
            M_k\mat{
            \alpha_{k+1,1}\\
            \vdots\\
            \alpha_{k+d,d}
            }=
            \mat{c^k_1\\
            \vdots\\
            c^k_d
            },
        \end{equation}
        where the entries of the matrix $M_k$ are constants and the entries of the vector $\mathbf{c}_k =(c^k_1,\dots, c^k_d)^\top$ are differential polynomials in the functions $u_\ell$ (all depending on $\chi$).
\end{lemma}

 \begin{remark}
    It is clear, from a geometric perspective, that the prescribed construction of the $\chi$-pentagram map yields a unique point $T_\e^\chi(\gamma)(x)$ for each point $\gamma(x)$. 
    In nondegenerate cases, one thus ought to be able to uniquely determine all $G_i$, and hence all coefficients $\alpha_{i,j}$. 
    Therefore, the matrix $M_k$ is invertible for nondegenerate configurations of $\chi$, while one may interpret the case when $M_k$ is noninvertible as a geometrically degenerate case.
\end{remark}

\begin{proof}[Proof of Lemma~\ref{lem:coef_computation}]
For $k=0$ we start by considering the first subspace $P_\e^1(x)=\big(\gamma(x+p_{1,0}\e),\dots, \gamma(x+p_{1,q_1}\e)\big)$, and expand the corresponding coplanarity condition \eqref{eq:chi_wedge_condition} in $\epsilon$. 
    For $1\leq j\leq d-q_1$, we see from the coefficient of $\Gamma \wedge \Gamma' \wedge \dots \wedge \Gamma^{(q_1)} \wedge \Gamma^{(q_1 + j)}$
    in the $\epsilon^{(0+1+\dots+q_1)+(q_1+j)}$ term that
    \begin{equation*} 
        0 = \sum \frac{p_{1,0}^{\beta_0} \dots p_{1, q_1}^{\beta_{q_1}}}{\beta_0! \dots \beta_{q_1}!} \cdot  \alpha_{\ell, \ell}\Gamma^{(\ell)} \wedge \Gamma^{(\beta_0)} \wedge \dots \wedge \Gamma^{(\beta_{q_1})},
    \end{equation*}
    where the summation runs over choices of $\{\ell, \beta_0, \dots, \beta_{q_1} \} = \{0, 1, \dots, q_1, q_1 + j \}$.
    We therefore obtain a linear relation of the form
    \begin{equation*}
        m^0_{j,1}\alpha_{1,1}+m^0_{j,2}\alpha_{2,2}+\dots+m^0_{j,q_1}\alpha_{q_1,q_1} + m^0_{j,q_1 + j}\alpha_{q_1+j,q_1 + j} = c^0_j,
    \end{equation*}
    where the coefficients $m^0_{j,\ell}$ and $c^0_j$ are constants in terms of $\chi$. (Here $c^0_j$ arises from the term with $\ell=0$.)
    Carrying this process out for each $P_\e^i$, we obtain $(d-q_1)+\dots+(d-q_r)=d$ such relations, from which we form the desired system \eqref{eq:M_system}.

    For $k > 0$, we modify the above process as follows: for $1\leq j\leq d-q_i$, now consider the $\Gamma\wedge \Gamma'\wedge \dots\wedge \Gamma^{(q_i)}\wedge \Gamma^{(q_i+j)}$ coefficient of the $\epsilon^{(0+1+\dots+q_i)+(q_i+j+k)}$ term.
    We obtain $d-q_i$ linear relations on $\alpha_{k+1,1},\alpha_{k+2,2},\dots,\alpha_{k+d,d}$, whose coefficients $m^k_{j,\ell}$ are again constants in terms of $\chi$, but now the terms $c^k_j$ are polynomials in terms of: 
    \begin{itemize}
        \item the parameters $\chi$;
        \item the functions $u_\ell$ and derivatives thereof, resulting from higher derivatives of $\Gamma$ being reduced through \eqref{eq:differentialeq}; and
        \item the previously-obtained coefficients, which by induction are also differential polynomials in the functions $u_\ell$.
    \end{itemize}
    Carrying this process out for each $P_\e^i$, we thus obtain a system \eqref{eq:M_system} of the desired form. This completes the proof
    of Lemma \ref{lem:coef_computation} and Proposition \ref{prop:chi_pent_coef_expansion}. 
\end{proof}

\begin{cor}\label{cor:coef_properties}
    The coefficients $\alpha_{i,i}$ are constants and $\alpha_{i+1,i}=0$ for $0 \leq i \leq d$. 
\end{cor}

\begin{proof}
    Since the entries of $M_0$ and $\mathbf{c}_0$ are constants (depending on $\chi$), it immediately follows that the coefficients $\alpha_{i,i}$ are themselves constant.
    Next, observe that $\mathbf{c}_1=\bf{0}$ (as a consequence of the fact that the differential equation \eqref{eq:differentialeq} does not involve $\Gamma^{(d)}$), and hence $\alpha_{i+1,i}=0$ for $1 \leq i \leq d$.
    Finally, the $\epsilon^1$ term of the normalization condition yields $\alpha_{1,0}=0$.
\end{proof}

We now proceed to the proof of the section's main result.

\begin{proof}[Proof of Theorem~\ref{thm:chi_pent_cont_limit}]
    By Lemma~\ref{lem:higher_coef_zero} and Corollary~\ref{cor:coef_properties}, $G_1 = \alpha_{1,1} \partial$ and $G_2 = \alpha_{2,0} + \alpha_{2,2} \partial^2$, with $\alpha_{1,1}$ and $\alpha_{2,2}$ constant.
    From the $\e^2$ term of the normalization $\det(\ge, \ge', \dots, \ge^{(d)}) \equiv 1$, we have
\begin{alignat*}{2}
	0&=|G_2\Gamma,\Gamma',\Gamma'',\dots,\Gamma^{(d)}|+\dots+&&|\Gamma,\Gamma',\Gamma'',\dots,G_2^{(d)}\Gamma|\\
	& &&+|\Gamma,\Gamma',\Gamma'',\dots,\Gamma^{(d-2)},G_1^{(d-1)}\Gamma,G_1^{(d)}\Gamma|\\
	&= (d+1)\cdot \alpha_{2,0}-2\alpha_{2,2}u_{d-1}+&&\hspace{-0.75em}\alpha_{1,1}^2u_{d-1},
\end{alignat*}
from which it follows that $\alpha_{2,0}=\frac{2}{d+1}(\alpha_{2,2}-\frac{1}{2}\alpha_{1,1}^2)u_{d-1}$. 
Thus
\begin{gather*}
    G_2=\alpha_{2,2}\l(\partial^2+\frac{2}{d+1} u_{d-1}  \r) -  \frac{\alpha_{1,1}^2}{d+1} u_{d-1},
\end{gather*}
as required.
In the case when $\chi$ is centralized, $\alpha_{1,1}=0$ and so
\begin{equation*}
	G_2=\alpha_{2,2}\left(\partial^2+\frac{2}{d+1} u_{d-1}\right) = \alpha_{2,2} \cdot Q_2,
\end{equation*}
which corresponds to the $(2,d+1)$-KdV flow (cf.~Proposition~\ref{prop:gamma_evolution_to_kdv}).
\end{proof}


\subsection{Centralizing the \texorpdfstring{$\chi$-pentagram}{χ-pentagram} map} 
\label{sub:centralizing_chi_map}

It follows from Equation~\eqref{eq:G_2_expansion} that the $(2,d+1)$-KdV equation only appears in Theorem~\ref{thm:chi_pent_cont_limit} for centralized configurations $\chi$. 
While a general condition for $\chi$ to be centralized is cumbersome (as it involves computing $M_0^{-1}\mathbf{c}_0$), we provide sufficient conditions for several important cases.

\subsubsection{Symmetric constructions}
Perhaps the simplest way to centralize the $\chi$-pentagram map is to choose $\chi$ invariant under negation, i.e., $\ge$ invariant under the symmetry $\e \to - \e$, which immediately implies that $G_i=0$ for odd~$i$. 
One can sometimes satisfy this condition by making a parameter shift (cf.~Remark~\ref{rem:short_diagonal_shift} for the short-diagonal map), but this approach is useful only for symmetrical constructions.

\subsubsection{Intersection of hyperplanes} 
\label{ssub:arbitrary_hyperplanes}
    Consider the intersection of $d$ arbitrary hyperplanes in $\RP^d$, i.e., the $\chi$-pentagram map with parameters
\begin{equation*}
    \chi = \{ \{p_{1,0}, \dots, p_{1,d-1} \}, \dots,  \{p_{d,0}, \dots, p_{d,d-1} \}\}.
\end{equation*}
Let $\sigma_{i,j}$ denote the $j$th elementary symmetric polynomial on $\{p_{i,0}, \dots ,p_{i,d-1}\}$.
The following discussion generalizes that of~\cite{Mari-Beffa:2013}, which focuses on the cases $d=3$ and $4$ (although unlike~\cite{Mari-Beffa:2013}, we do not discuss integrability).

\begin{prop}[cf.~Proposition 4.3 of \cite{Mari-Beffa:2013}] \label{prop:plane_sym_poly_condition}
    Let $\chi$ describe the intersection of $d$ hyperplanes. Then the following conditions are equivalent:
    \begin{enumerate}
        \item The top elementary polynomials $\sigma_{i,d}$ $(1 \leq i \leq d)$ associated with all the hyperplanes coincide (denote their common value by $\sigma_d$).
        \item The coefficients $\alpha_{j,j} = 0$ for $j = 1, \dots, d-1$.
    \end{enumerate}
    Furthermore, if the above conditions hold, then $\alpha_{d,d} = (-1)^{d+1} \sigma_d/d!$.
\end{prop}
\begin{proof}
As detailed in Lemma~\ref{lem:coef_computation}, we have a system
\begin{equation}\label{eq:d_planes_system}
        M_0\mat{
        \alpha_{1,1}\\
        \vdots\\
        \alpha_{d,d}
        }= \mathbf{c}_0, 
    \end{equation}
the $i$th row of which is the linear relation coming from the coplanarity condition
\begin{equation*}
    \det(\ge(x), \Gamma(x + p_{i,0} \epsilon), \dots, \Gamma(x + p_{i,d-1} \epsilon)) = 0.
\end{equation*}
Rescaling the $i$th row by $-\prod_{k=0}^d k! / \prod_{k<\ell}(p_{i,\ell}-p_{i,k})$, we compute 
\begin{equation}\label{eq:M_0_entries}
    (M_0)_{i,j} = (-1)^{j+1}  j! \cdot  \sigma_{i, d-j} \quad \text{and } \quad (\mathbf{c}_0)_i = \sigma_{i,d}.
\end{equation}

Since all entries of the last column of $M_0$ are the constant $(-1)^{d+1} d!$, it follows that the $j$th row of ${M_0^{-1}}$ sums to $0$ for $1 \leq j \leq d-1$ and to $(-1)^{d+1}/d!$ for $j=d$. If $\sigma_{i,d} = \sigma_d$ for each $i$, we obtain the desired result by computing
\begin{equation*}
        \mat{
        \alpha_{1,1}\\
        \vdots\\
        \alpha_{d,d}
        }= M_0^{-1} \mat{ \sigma_d \\ \vdots \\ \sigma_d}.
    \end{equation*} 
    Conversely, if $\alpha_{1,1}=\dots=\alpha_{d-1,d-1} = 0$, then \eqref{eq:d_planes_system} implies that $\sigma_{i,d} = (-1)^{d+1} d! \cdot \alpha_{d,d}$ for each~$i$.
\end{proof} 

This is a fairly restrictive condition to impose on $\chi$.
One can obtain more explicit conditions for configurations of hyperplanes where there are additional combinatorial relationships between the $\sigma_{i,j}$, as in the following examples.

\subsubsection{Evenly spaced hyperplanes} 
\label{ssub:evenly_spaced}
Consider the intersection of $d$ evenly spaced hyperplanes in $\RP^d$, i.e., the $\chi$-pentagram map $T^\chi_\epsilon$ with 
\begin{equation*}
    \chi = \{ \{p_0 + ir, \dots, p_{d-1}+ir \} : i = 0, \dots, d-1 \}
\end{equation*}
for some $r \in \R$. 
When the parameters $p_i$ and $r$ are all integers, $T^\chi_\epsilon$ is a continuous analogue of the generalized pentagram map $T_{I, J}$ from~\cite{Khesin:2016}, with $J = (r,\dots, r)$. In particular, this includes the case of the short-diagonal map.
Using~\eqref{eq:M_0_entries} we can compute
\begin{equation*}
    \alpha_{1,1} = \frac{1}{d} \left(\sigma_1 + \binom{d}{2} r \right),
\end{equation*}
where $\sigma_1 = \sum_{i=0}^{d-1} p_i$. Thus $\chi$ is centralized if and only if $\sigma_1 = - \binom{d}{2}r$ (cf.~the condition in~\cite{Khesin:2016} for the envelope construction of the continuous limit). 
For instance, the parameter shifts specified in Equation~\eqref{eq:shifted_short_diag} for the short-diagonal map ensure that this condition is satisfied.
    


\subsubsection{Dual dented map} 
\label{ssub:dual_dented_map}

Let $s \in \{ 1, \dots, d-1 \}$. As described in Example~\ref{ex:dual_dented}, the continuous analogue of the dual dented map $\widehat{T}_s$ is given by the map $T_\e^{\widehat{\chi}_s}$ with parameters
\begin{equation*}
    \widehat{\chi}_s \coloneqq \{ \{i, i+1, \dots, i + d-1  \} : 0 \leq i \leq d,\,  i \neq d-s \}.
\end{equation*}
We seek to centralize the dual dented map by shifting $\widehat{\chi}_s$ by some $\delta \in \R$ (i.e.~replacing each $p_{i,j}$ by $p_{i,j} + \delta$); this is equivalent to considering the shifted curve
\begin{equation*}
    \tilde{\gamma}_\epsilon(x) = T_\e^{\widehat{\chi}_s} (\gamma) (x + \delta \epsilon)
\end{equation*}
(cf.~Equation~\eqref{eq:shifted_short_diag} for the short-diagonal map).

\begin{prop}\label{prop:dual_dented_limit}
    The dual dented map becomes centralized after shifting $\widehat{\chi}_s$ by $1-d-\frac{s}{d}$, and hence its continuous limit is then the $(2,d+1)$-KdV equation.
\end{prop}

The proof of this proposition follows from a lengthy but straightforward combinatorial computation.

\begin{example}
Consider the case $d=3, s=1$, shown in Figure~\ref{fig:dual_dented}. Note that the shifted point $\tilde{\gamma}_\epsilon (x + \frac{7}{3}) = T_\epsilon^{\widehat{\chi}_1}(\gamma)(x)$ is centred near $\gamma (x + \frac{7}{3} )$, which agrees with the shift by $1-3-\frac 13=-\frac 73$ given in 
Proposition \ref{prop:dual_dented_limit}.
\end{example}



\section{Continuous limit of the Lax form}\label{sec:lax_forms} 

The short-diagonal pentagram map is known to have a Lax representation in both the discrete and continuous case.
We will apply quantum calculus to demonstrate that the Lax matrices in the discrete case tend to those in the continuous case, and show that the continuous Lax equation can be obtained as a limit of the discrete one.
First, we recall the relevant notions.

\subsection{Lax representations of the short-diagonal pentagram map} 
\label{sub:lax_review}

In the continuous case, a \emph{Lax equation} is a differential equation of the form $\frac{d}{dt} L = [V, L]$, where $L$ and $V$ are time-dependent linear (i.e.~differential) operators. 
Suppose $L$ is a matrix first-order differential operator $L = \frac{d}{dx} - U$.
Then the Lax equation assumes the form of the \emph{zero-curvature equation}
\begin{equation} \label{eq:zero-curv} 
    \frac{d}{dt} U = [V, U] + \frac{d}{dx} V.
\end{equation}
This is the compatibility condition of the system of PDEs
\begin{equation*} \label{eq:matrix-pde-system}
    \begin{cases}
        \frac{d}{dx} \Psi = U \Psi \\
        \frac{d}{dt} \Psi = V \Psi,
    \end{cases}
\end{equation*} 
where $U=U(x,t),\, V=V(x,t)$ are $(d+1)\times(d+1)$ matrices and $\Psi = (\psi, \psi', \dots, \psi^{(d)})^\top$.

\begin{example}
    As mentioned in Remark~\ref{rem:kdv-lax}, the $(2,d+1)$-KdV equation arising from the continuous limit of the short-diagonal map can be written in the Lax form. In this case we have $L = \partial^{d+1} + u_{d-1} \partial^{d-1} + \dots + u_0$ and $V = Q_2 \coloneqq (L^{2/d+1})_+$, satisfying $dL/dt = [V,L]$. 
To simplify the calculations below, we will set $V = c \cdot Q_2$, i.e.~we consider the evolution
\begin{equation*}
    \frac{d}{dt} \Gamma = c \, Q_2 \Gamma,
\end{equation*}
where $c = \tilde{C}_d$ is the constant from Theorem~\ref{thm:directctslimit}, which amounts to rescaling the time parameter.
In the matrix formulation, we have
\begin{equation} \label{eq:U-matrix}
    U = \left( \begin{array}{c|cccc}
    0 & \\
    \vdots & \multicolumn{4}{c}{I} \\
    0 & \\
    \hline
    -u_0 & -u_1 &  \cdots & - u_{d-1} & 0
    \end{array} \right),
\end{equation}
while $V$ is the unique matrix such that
\begin{equation} \label{eq:V-matrix}
    V \begin{pmatrix}
        \Gamma \\
        \Gamma' \\
        \vdots \\
        \Gamma^{(d)}
    \end{pmatrix} = \frac{d}{dt} \begin{pmatrix}
        \Gamma \\
        \Gamma' \\
        \vdots \\
        \Gamma^{(d)}
    \end{pmatrix} = 
    c \begin{pmatrix}
        Q_2 \Gamma \\
        (Q_2\Gamma)' \\
        \vdots \\
        (Q_2 \Gamma)^{(d)}
    \end{pmatrix}
\end{equation}
(recall that $\Gamma(x), \Gamma'(x), \dots, \Gamma^{(d)}(x)$ form a basis for each $x$).
\end{example}

A \emph{discrete Lax equation} (or \emph{discrete zero-curvature equation}) with spectral parameter is an equation of the form
\begin{equation} \label{eq:discrete-lax}
    L_{i,t+1}(z) = P_{i+1,t}(z) L_{i,t}(z) P_{i,t}^{-1}(z),
\end{equation}
where $i,t \in \Z_{\geq 0}$ and $z \in \mathbb{C}$ is the spectral parameter. Analogously to the continuous case, this is the compatibility equation ensuring a solution $\psi_{i,t}(z)$ to the overdetermined system
\begin{equation*}
    \begin{cases}
         L_{i,t}(z) \psi_{i,t}(z) = \psi_{i+1, t}(z) \\
         P_{i,t}(z) \psi_{i,t}(z) = \psi_{i, t+1}(z).
    \end{cases} 
\end{equation*}

\begin{examplecont}[continued]
    A discrete Lax form for the short-diagonal pentagram map is described in~\cite{Khesin:2013}, where it is used to establish integrability in the algebraic-geometric sense. 
The construction of the Lax matrix $L_{i,t}(z)$ uses scaling invariance of the pentagram map, which was proved in~\cite{Khesin:2013, Mari-Beffa:2015}.

In this situation, the index $i$ corresponds to vertices of a twisted $n$-gon, and each increment of the index $t \to t+1$ corresponds to an iteration of the pentagram map.
Fix a time $t\in \Z_{\geq 0}$ and a twisted $n$-gon with vertices $(v_i)$.
It will be convenient for our purposes to replace the matrices $P_{i,t}$ and $L_{i,t}$ from~\cite{Khesin:2013} with their transpose-inverses (still satisfying \eqref{eq:discrete-lax}). Namely, we consider
\begin{equation*}
    L_{i,t}(z) = \left( \begin{array}{c|ccc}
    0 & \\
    \vdots & & \Lambda(z)\\
    0 & \\
    \hline
    (-1)^d & a_{i,1}&  \cdots& a_{i,d}
    \end{array} \right)
\end{equation*}
where 
\begin{equation*}
    \Lambda(z) = \begin{cases}
        \text{diag}(z,1,z,1, \dots, z) & \text{if $d$ is odd,}  \\
        \text{diag}(1,z,1,z,\dots,1,z) & \text{if $d$ is even,} 
    \end{cases}
\end{equation*}
and the $a_{i,j}$ are the coordinates associated to the canonical lift of the $v_i$ (see Equation~\eqref{eq:differenceeq}). 
An explicit formula for $P_{i,t}(z)$ is given in~\cite{Soloviev:2013} for the 2D case and in~\cite{Khesin:2013} for the 3D case, but is not known in general.
However, we will be able to circumvent this issue by using our direct construction of the continuous limit from Section~\ref{sub:limit_of_short_diag_pentagram_map}.
\end{examplecont}

Our main result in this section is the following analysis of the continuous limit:
\begin{theorem}\label{thm:lax_limit}
    In the continuous limit as $n \to \infty$, the space of Lax matrices $L_{i,t}(z)$ at $z=1$ tends to the space of matrices $U$ corresponding to differential operators of the form $\partial^{d+1} + u_{d-1} \partial^{d-1} + \dots + u_0$.
    Furthermore, in the limit, the dynamics of the discrete Lax matrices become the dynamics of the differential operators, as described by the $(2, d+1)$-KdV zero-curvature equation.
\end{theorem}

Analogously to Section~\ref{sec:ctslimitcalc}, we will study the continuous limit explicitly by associating a Lax matrix $\tilde{L}_{i,t}$ to the discretization of a lifted curve. 
The two parts of the theorem are stated more technically and proved as Propositions \ref{prop:lax-matrix-expansion} and \ref{prop:limit-of-lax-eq} below.

\subsection{Kinematics: Limit of the Lax matrices} 
\label{sub:limit_of_the_lax_matrix}

Following the same discretization procedure as in Section~\ref{sec:ctslimitcalc}, we fix a nondegenerate curve $\gamma$ in $\RP^d$, let $v_j \coloneqq \gamma(x+j \epsilon)$ for fixed $x$ and small $\epsilon$, and lift $\gamma$ to $\Gamma$ in $\R^{d+1}$ such that $\det (\Gamma, \Gamma', \dots, \Gamma^{(d)}) \equiv 1$. 
Then $\Gamma$ satisfies a differential equation $L \Gamma = 0$ where $L = \partial^{d+1} + u_{d-1} \partial^{d-1} +\dots +  u_0$, and 
the points $\tilde{V}_j \coloneqq \Gamma(x+ j \epsilon)$ satisfy difference equations \eqref{eq:diffeq-v1} and \eqref{eq:diffeq-v2} in coordinates $\tilde{a}_k$ and $A_k$ respectively.

We introduce a discrete Lax matrix associated to the points $\Gamma(x), \Gamma(x+\epsilon), \dots, \Gamma(x + d \epsilon)$:
\begin{equation*}
    \tilde{L}_{0,0}(z) \coloneqq \left( \begin{array}{c|ccc}
    0 & \\
    \vdots & & \Lambda(z)\\
    0 & \\
    \hline
    \tilde{a}_0 & \tilde{a}_1 &  \cdots& \tilde{a}_d
    \end{array} \right), \qquad \tilde{L}_{0,0} \coloneqq \tilde{L}_{0,0}(1) = \left( \begin{array}{c|ccc}
    0 & \\
    \vdots & & I \\
    0 & \\
    \hline
    \tilde{a}_0 & \tilde{a}_1 &  \cdots& \tilde{a}_d
    \end{array} \right).
\end{equation*}
(Note that since the discretization is local near $\gamma(x)$ and we are only considering the initial polygon for now, we may assume that $i=t=0$.) 

\begin{remark}
    Whereas the earlier Lax matrix $L_{0,0}(1)$ was associated to the canonical lift of the vertices $v_i$, the matrix $\tilde{L}_{0,0}$ is associated to the canonical lift of the curve $\gamma \to \Gamma$, and as such can be regarded as an approximation of the former matrix  (cf.~Remark~\ref{rem:vertices_vs_curve_lift}). 
\end{remark}

In order to study the expansion of $\tilde{L}_{0,0}$ as $\epsilon \to 0$, we define a change of basis matrix $D_\epsilon$ by
\begin{equation*}
    D_\epsilon \begin{pmatrix}
        \Gamma(x)\\
        \Gamma(x + \epsilon) \\
        \vdots \\
        \Gamma(x+d \epsilon)
    \end{pmatrix} = \begin{pmatrix}
        \Gamma(x) \\
        \Delta_\epsilon \Gamma(x)/\epsilon \\
        \vdots  \\
        \Delta_\epsilon^d \Gamma(x ) / \epsilon^d
    \end{pmatrix},
\end{equation*}
where $\Delta_\epsilon$ is the difference operator introduced in Section~\ref{sec:ctslimitcalc} 
(cf.~the rewriting of shift operators as difference operators for Proposition~\ref{prop:limit-of-difference-eq}).

\begin{prop}\label{prop:lax-matrix-expansion}
    Expanded as a series in $\epsilon$, we have
    \begin{equation}\label{eq:L-expansion}
        D_\epsilon \tilde{L}_{0,0} D_\epsilon^{-1} = I + \epsilon U + \O(\epsilon^2),
    \end{equation}
    where $U = U(x,0)$. 
\end{prop}
\begin{proof}
First note that
\begin{equation*}
    \tilde{L}_{0,0} \begin{pmatrix}
    	\Gamma(x)\\
    	\vdots \\
    	\Gamma(x+d \epsilon)
    \end{pmatrix} = \begin{pmatrix}
    	\Gamma(x+\epsilon) \\
    	\vdots  \\
    	\Gamma(x + (d+1) \epsilon)
    \end{pmatrix}
\end{equation*}
by the difference equation \eqref{eq:diffeq-v1}.
Then one computes 
\begin{equation*} \label{eq:discretized-derivative-matrix}
    (D_\epsilon \tilde{L}_{0,0} D_\epsilon^{-1} - I) \begin{pmatrix}
    	\Gamma(x) \\
    	\Delta_\epsilon \Gamma(x)/\epsilon \\
    	\vdots  \\
    	\Delta_\epsilon^d \Gamma(x ) / \epsilon^d
    \end{pmatrix} = \begin{pmatrix}
    	\Delta_\epsilon \Gamma(x) \\
    	\Delta_\epsilon^2 \Gamma(x)/\epsilon \\
    	\vdots  \\
    	\Delta_\epsilon^{d+1} \Gamma(x) / \epsilon^d
    \end{pmatrix},
\end{equation*}
from which it follows that
\begin{equation*}
    \epsilon^{-1} (D_\epsilon \tilde{L}_{0,0} D_\epsilon^{-1} - I) = \left( \begin{array}{c|cccc}
    0 & \\
    \vdots & \multicolumn{4}{c}{I} \\
    0 & \\
    \hline
    -A_0/\epsilon^{d+1} & -A_1/\epsilon^d &  \cdots& - A_{d-1}/\epsilon^2 & -A_d/\epsilon
    \end{array} \right)
\end{equation*}
by the difference equation \eqref{eq:diffeq-v2}.

Using the expansions of the $A_i$ from Lemma~\ref{lem:A_i expansions}, we see that
\begin{equation*}
    \lim_{\epsilon \to 0} \epsilon^{-1} (D_\epsilon \tilde{L}_{0,0} D_\epsilon^{-1} - I) = \left( \begin{array}{c|cccc}
    0 & \\
    \vdots & \multicolumn{4}{c}{I} \\
    0 & \\
    \hline
    -u_0 & -u_1 &  \cdots & - u_{d-1} & 0
    \end{array} \right) = U,
\end{equation*} 
as required. 
\end{proof}

\begin{remark}\label{rmk:L-matrix-refinement}
    Later on we will need a slightly refined description of the higher-order terms of~\eqref{eq:L-expansion}.
    It follows from Remark~\ref{rmk:A_i refinement} (which describes the higher-order terms of the $A_i$) that all of the coefficients in the power series expansion for $D_\epsilon \tilde{L}_{0,0} D_\epsilon^{-1}$ can be written as polynomial functions of the $u_j$ and their derivatives. 
\end{remark}

\subsection{Dynamics: Limit of the Lax equation}
\label{sub:limit_of_the_lax_equation}

So far, we have defined a matrix $\tilde{L}_{0,0}$ associated to the discretization of a curve $\Gamma$ at time $t=0$. 
We extend this to define $\tilde{L}_{i,0}$ for any $i$, by replacing $x$ with $x+ i \epsilon$ in the discretization procedure.
Next, define $\tilde{P}_{i,0}$ to be the unique matrix such that
\begin{equation*}
    \tilde{P}_{i,0} \begin{pmatrix}
        \Gamma(x+i \epsilon)\\
        \vdots \\
        \Gamma(x+(i+d) \epsilon)
    \end{pmatrix} = \begin{pmatrix}
        \tilde{\Gamma}_\epsilon(x+i \epsilon)\\
        \vdots \\
        \tilde{\Gamma}_\epsilon(x+(i+d) \epsilon)
    \end{pmatrix},
\end{equation*}
where $\tilde{\Gamma}_\epsilon$ is the curve obtained via the pentagram map (see Section~\ref{sub:limit_of_short_diag_pentagram_map}). 
Finally, we define $\tilde{L}_{i,t}$ and $\tilde{P}_{i,t}$ for all $i$ and $t$ by inductively replacing $\Gamma$ with $\tilde{\Gamma}_{\epsilon}$ in the above definitions when we increment $t \to t+1$. 
By construction, the matrices $\tilde{L}_{i,t}$ and $\tilde{P}_{i,t}$ satisfy the discrete Lax equation
\begin{equation} \label{eq:discrete-lax-tilde}
    \tilde{L}_{i,t+1} = \tilde{P}_{i+1,t} \tilde{L}_{i,t} \tilde{P}_{i,t}^{-1},
\end{equation}
and can be viewed as approximations of the matrices $L_{i,t}(1)$ and $P_{i,t}(1)$ associated to twisted $n$-gons. 

\begin{prop}\label{prop:limit-of-lax-eq}
    The zero-curvature equation \eqref{eq:zero-curv} for the continuous limit of the pentagram map can be obtained as a limit of the discrete Lax equation \eqref{eq:discrete-lax-tilde}. 
    More precisely, by rewriting \eqref{eq:discrete-lax-tilde} for $i=t=0$ as
    \begin{equation*}
        \frac{D_\epsilon \tilde{L}_{0,1} D_\epsilon^{-1} - D_\epsilon \tilde{L}_{0,0} D_\epsilon^{-1}}{\epsilon^3} = \frac{(D_\epsilon \tilde{P}_{1,0} D_\epsilon^{-1})(D_\epsilon \tilde{L}_{0,0} D_\epsilon^{-1})(D_\epsilon \tilde{P}_{0,0} D_\epsilon^{-1})^{-1} - D_\epsilon \tilde{L}_{0,0} D_\epsilon^{-1}}{\epsilon^3},
    \end{equation*}
    we obtain the corresponding KdV zero-curvature representation
    \begin{equation*}
        \frac{d}{dt} U = [V, U] + \frac{d}{dx} V
    \end{equation*}
    for $U$ and $V$ as in \eqref{eq:U-matrix} and \eqref{eq:V-matrix}, as $\epsilon \to 0$.
\end{prop}

This statement only involves $t=0$ and $t=1$, i.e., a single iteration of the pentagram map (parametrized by $\epsilon$).
Recall that the curve $\tilde{\Gamma}_\epsilon$ at $t=1$ satisfies a differential equation ${L}_\epsilon \tilde{\Gamma}_\epsilon = 0$, where $L_\epsilon = \partial^{d+1} + u_{d-1, \epsilon} \partial^{d-1} + \dots + u_{0, \epsilon}$. 
By expanding $u_{i,\epsilon} = u_i + \epsilon^2 w_i + \O(\epsilon^4)$, we (by definition) obtain the evolution $du_i/dt = w_i$, regarding $\epsilon^2$ as time.
Let
 \begin{equation*}
    U_\epsilon \coloneqq \left( \begin{array}{c|cccc}
    0 & \\
    \vdots & \multicolumn{4}{c}{I} \\
    0 & \\
    \hline
    -u_{0,\epsilon} & -u_{1, \epsilon} &  \cdots & - u_{d-1, \epsilon} & 0
    \end{array} \right). 
\end{equation*}
Our first step is to relate
\begin{equation*}
    \frac{d}{dt}U = \lim_{\epsilon \to 0} \frac{U_\epsilon - U}{\epsilon^2}
\end{equation*}
to the matrices $\tilde{L}_{0,0}$ and $\tilde{L}_{0,1}$.

\begin{lemma} \label{lem:limit-of-lax-1}
    $D_\epsilon \tilde{L}_{0,1} D_\epsilon^{-1} - D_\epsilon \tilde{L}_{0,0} D_\epsilon^{-1} = \epsilon(U_\epsilon - U) + \O(\epsilon^4)$.
\end{lemma}
\begin{proof}
    As mentioned in Remark~\ref{rmk:L-matrix-refinement}, one can write 
\begin{equation*}
    D_\epsilon \tilde{L}_{0,0} D_\epsilon^{-1} = I + \epsilon U + \epsilon^2 F_2(u_0, \dots, u_{d-1}) + \epsilon^3 F_3(u_0, \dots, u_{d-1})  + \O(\epsilon^4),
\end{equation*}
where $F_2$ and $F_3$ are polynomial functions of the $u_j$ and their derivatives.
Replacing $\Gamma$ with $\tilde{\Gamma}_\epsilon$ (and hence $u_j$ with $u_{j, \epsilon}$) in the definition of $\tilde{L}_{0,0}$, we get
\begin{align*}
    D_\epsilon \tilde{L}_{0,1} D_\epsilon^{-1} &= I + \epsilon U_\epsilon + \epsilon^2 F_2(u_{0, \epsilon}, \dots, u_{d-1, \epsilon}) + \epsilon^3 F_3(u_{0, \epsilon}, \dots, u_{d-1, \epsilon})  + \O(\epsilon^4) \\
    &= I + \epsilon U_\epsilon + \epsilon^2 F_2(u_0, \dots, u_{d-1}) + \epsilon^3 F_3(u_0, \dots, u_{d-1})  + \O(\epsilon^4)
\end{align*}
where the second equality uses that $u_{j, \epsilon} = u_j + \O(\epsilon^2)$.
\end{proof}

Next, we study the expansion of the terms of the Lax equation involving $\tilde{P}_{0,0}$ and $\tilde{P}_{1,0}$.

\begin{lemma} \label{lem:limit-of-lax-2}
    Expanded as series in $\epsilon$, we have
    \begin{equation} \label{eq:P-expansion-0}
        D_\epsilon \tilde{P}_{0,0} D_\epsilon^{-1} = I + \epsilon^2 V + \O(\epsilon^4),
    \end{equation}
    \begin{equation} \label{eq:P-expansion-1}
        (D_\epsilon \tilde{P}_{0,0} D_\epsilon^{-1})^{-1} = I - \epsilon^2 V + \O(\epsilon^4),
    \end{equation}
    and
    \begin{equation}\label{eq:P-expansion-2}
        D_\epsilon \tilde{P}_{1,0} D_\epsilon^{-1} = I + \epsilon^2 V + \epsilon^3 \frac{d}{dx}V + \O(\epsilon^4),
    \end{equation}
    where $V = V(x,0)$.
\end{lemma}
\begin{remark}
    This key lemma allows one to identify the limit of the discrete Lax equation as a familiar KdV-type expression in any dimension $d$, even though an explicit form of the matrices $P_{i,t}$ remains obscure beyond $d=2$ and $3$.
\end{remark}
\begin{proof}
The first statement follows from the expansion $\tilde{\Gamma}_\epsilon = \Gamma + \epsilon^2 c \,Q_2 \Gamma + \O(\epsilon^4)$ given in Theorem~\ref{thm:directctslimit}, where $Q_2 = \partial^2 + \frac{2}{d+1} u_{d-1}$ and $c = \tilde{C}_d$. 
Indeed, one computes
\begin{equation*} 
    \epsilon^{-2} (D_\epsilon \tilde{P}_{0,0}D_\epsilon^{-1} - I) \begin{pmatrix}
        \Gamma(x)\\
        \Delta_\epsilon \Gamma(x)/\epsilon \\
        \vdots \\
        \Delta_\epsilon^d \Gamma(x) /\epsilon^d
    \end{pmatrix} = c\begin{pmatrix}
        Q_2 \Gamma(x)\\
        \Delta_\epsilon Q_2 \Gamma (x)/\epsilon \\
        \vdots \\
        \Delta_\epsilon^d Q_2 \Gamma(x)/\epsilon^d
        \end{pmatrix} + \O(\epsilon^2),
\end{equation*}
and therefore
\begin{equation*}
    \lim_{\epsilon\to 0 } \epsilon^{-2} (D_\epsilon \tilde{P}_{0,0}D_\epsilon^{-1} - I) \begin{pmatrix}
        \Gamma(x)\\
        \Gamma'(x) \\
        \vdots \\
        \Gamma^{(d)}(x)
    \end{pmatrix} = c \begin{pmatrix}
        Q_2 \Gamma(x)\\
        (Q_2 \Gamma)' (x) \\
        \vdots \\
        (Q_2 \Gamma)^{(d)}(x)
        \end{pmatrix}.
\end{equation*}
But $V$ was defined as the unique matrix satisfying the above equation, so we must have $V = \lim_{\epsilon\to 0 } \epsilon^{-2} (D_\epsilon \tilde{P}_{0,0}D_\epsilon^{-1} - I)$.
This gives the desired expansion \eqref{eq:P-expansion-0}, and inverting that power series yields \eqref{eq:P-expansion-1}. 
Finally, by replacing $x$ with $x + \epsilon$ in \eqref{eq:P-expansion-0}, we get 
\begin{equation*}
    D_\epsilon \tilde{P}_{1,0} D_\epsilon^{-1} = I + \epsilon^2 V(x+\epsilon,0) + \O(\epsilon^4),
\end{equation*}
and we obtain \eqref{eq:P-expansion-2} by expanding $V(x+\epsilon, 0)$ in $\epsilon$.
\end{proof}

It is now straightforward to calculate the limits in Proposition~\ref{prop:limit-of-lax-eq}.

\begin{proof}[Proof of Proposition~\ref{prop:limit-of-lax-eq}] 
    For the first limit, using Lemma~\ref{lem:limit-of-lax-1}, we have
    \begin{align*}
        \lim_{\epsilon\to 0}\frac{D_\epsilon \tilde{L}_{0,1} D_\epsilon^{-1} - D_\epsilon \tilde{L}_{0,0} D_\epsilon^{-1}}{\epsilon^3}  &= \lim_{\epsilon \to 0} \frac{\epsilon(U_\epsilon - U) + \O(\epsilon^4)}{\epsilon^3} = \lim_{\epsilon \to 0} \frac{U_\epsilon - U}{\epsilon^2} = \frac{d}{dt} U.
    \end{align*} 
    For the second limit, using the expansions \eqref{eq:P-expansion-1} and \eqref{eq:P-expansion-2} from Lemma \ref{lem:limit-of-lax-2} and writing $D_\epsilon \tilde{L}_{0,0} D_\epsilon^{-1} = I + \epsilon U +  R$ for some remainder term $R = \O(\epsilon^2)$, we have
    \begin{align*}
        &\, \lim_{\epsilon \to 0} \frac{(D_\epsilon \tilde{P}_{1,0} D_\epsilon^{-1})(D_\epsilon \tilde{L}_{0,0} D_\epsilon^{-1})(D_\epsilon \tilde{P}_{0,0} D_\epsilon^{-1})^{-1} - D_\epsilon \tilde{L}_{0,0} D_\epsilon^{-1}}{\epsilon^3} \\
        =&\,  \lim_{\epsilon \to 0} \frac{(I + \epsilon^2 V + \epsilon^3 \frac{d}{dx} V + \O(\epsilon^4))(I + \epsilon U + R)(I - \epsilon^2 V + \O(\epsilon^4)) - (I + \epsilon U + R)}{\epsilon^3} \\
        =&\,  \lim_{\epsilon \to 0}  \frac{\epsilon^3 (VU - UV + \frac{d}{dx} V) + \O(\epsilon^4)}{\epsilon^3} = [V, U] + \frac{d}{dx} V.
    \end{align*}
    Thus we have recovered the continuous zero-curvature equation as a limit of the discrete Lax equation.
\end{proof}


\section{Realization of the KdV-type equations}\label{sec:kdv_realization}

Expanding upon the results of~\cite{Mari-Beffa:2013}, we obtain specific instances of $\chi$ for which the corresponding $\chi$-pentagram map evolution yields various KdV flows through its continuous limit. Furthermore, we provide a heuristic evidence that not all KdV equations may be realized through pentagram-type maps, thus partially answering questions posed in~\cite{Khesin:2013} and~\cite{Mari-Beffa:2013}.

\subsection{Low KdV equations}
By computing the $\alpha_{i,j}$ coefficients via the procedure detailed in Lemma~\ref{lem:coef_computation} and setting these equal to the coefficients of $(L^{m/d+1})_+$, we obtain a system in terms of the parameters $p_{i,j}$ of $\chi$. Solutions of this system, whenever they exist, yield a configuration of~$\chi$ for which the continuous limit of $T_\e^\chi$ corresponds to the desired $(m,d+1)$-KdV equation. This computation for low values of $m,d$ is done in~\cite{Mari-Beffa:2013}, and we exemplify here a more general approach by realizing the $(3,4)$-KdV equation.

To construct the  $(3,4)$-KdV equation we consider the $\chi$-pentagram map with
\begin{equation*}
    \chi=\l\{\{p_{1,0},p_{1,1},p_{1,2} \},\{p_{2,0},p_{2,1},p_{2,2} \},\{p_{3,0},p_{3,1},p_{3,2} \}\r\},
\end{equation*}
corresponding to the intersection of three $2$-dimensional planes in $\RP^3$. We describe all possible $\chi$ of this form which yield the $(3,4)$-KdV equations. Our approach may be regarded as an alternative to that in Section 4.1 of \cite{Mari-Beffa:2013} in somewhat more invariant terms. For simplicity, let $\sigma_{i,j}$ denote the $j$th elementary symmetric polynomial on $\{p_{i,0},p_{i,1},p_{i,2}\}$. 

\begin{prop}\label{prop:conditions_for_(3,4)}
    With $\chi$ as above, the $\e$-expansion of $\ge$ will be of the form 
    \begin{equation*}
        \ge=\Gamma+\e^3c \cdot \big(L^{3/4}\big)_{+}\Gamma+\dots
    \end{equation*}
    for some constant $c\neq 0$ if and only if both $\sigma_{1,3}=\sigma_{2,3}=\sigma_{3,3}\neq 0$ and one polynomial constraint in $\sigma_{i,j}$, described below and corresponding to $\alpha_{3,1}=\frac{3}{4}u_2\alpha_{3,3}$, are satisfied.
\end{prop}

\begin{proof}
The expansions \eqref{eq:G_1_expansion} and \eqref{eq:G_2_expansion} show that $G_1$ and $G_2$ vanish if and only if $\alpha_{1,1} = \alpha_{2,2} = 0$.
Therefore, it follows from Proposition \ref{prop:plane_sym_poly_condition} that $G_1=G_2=0$, $G_3\neq 0$ if and only if $\sigma_{1,3}=\sigma_{2,3}=\sigma_{3,3}\neq 0$.
Next, as described in Lemma \ref{lem:coef_computation}, we may compute $\alpha_{3,1}$ by evaluating $M_2^{-1}\mathbf{c}_2$, from which we rewrite $\alpha_{3,1}=\frac{3}{4}u_2\alpha_{3,3}$ (coming from $(L^{3/4})_+$) in terms of $\sigma_{i,j}$. After satisfying this final constraint, one uses 
the normalization on $\ge$ to solve for $\alpha_{3,0}$ and find
$G_3=\alpha_{3,3}\big(L^{3/4}\big)_+$
as desired. 
We have thus obtained polynomial conditions, namely $\sigma_{1,3}=\sigma_{2,3}=\sigma_{3,3}\neq 0$ and the rewriting of $\alpha_{3,1}=\frac{3}{4}u_2\alpha_{3,3}$ in terms of $\sigma_{i,j}$, which will be met if and only if the $\e$-expansion of $\ge$ is of the prescribed form. 
\end{proof}

\begin{example}
Theorem 4.6 of~\cite{Mari-Beffa:2013} provides a sufficient condition for 
$\chi$ to give rise to the $(3,4)$-KdV equation in the continuous limit. 
Namely, the configuration $\chi$ could be of the form $\chi=\{\{-c,a,b\},\{c,-a,b\},\{c,-1,ab\}\}$ satisfying the relation
$c-1+a(b-1)=-5(b-c)/4$, such as e.g.~an integer solution $a=-2, b=3, c=-5$. 

The conditions described in Proposition~\ref{prop:conditions_for_(3,4)} generalize the condition from ~\cite{Mari-Beffa:2013}, as may be checked through a tedious computation. 
For instance, taking
\begin{equation*}
    \chi=\l\{\{-1, \tfrac{3}{2}, 4\},\{\tfrac{6}{5}, 10, -\tfrac{1}{2}\},\{1, -r, \tfrac{6}{r}\}\r\},
\end{equation*}
where $r$ is any of the four real roots of the polynomial
\begin{equation*}
 R(x)=   2480 x^4+33006 x^3+72121 x^2-198036 x+89280,
\end{equation*}
satisfies the conditions of Proposition \ref{prop:conditions_for_(3,4)} and hence yields the $(3,4)$-KdV equation, yet lies outside the condition described in~\cite{Mari-Beffa:2013}. 
A geometric interpretation of constraints in Proposition~\ref{prop:conditions_for_(3,4)} is not immediately apparent, while the existence of integer solutions satisfying them but lying outside of the sufficient condition  in \cite{Mari-Beffa:2013} is an open question.
\end{example}

\begin{remark}
The $(3,5)$-KdV equation was realized in~\cite{Mari-Beffa:2013}. Similarly, one would expect that the above approach may generate examples of a $\chi$-pentagram map with the intersection of 4 hyperplanes in $\RP^4$ for which the continuous limit is the $(4,5)$-KdV equation. However, finding such a configuration (in particular, one with integer coefficients) remains an open problem.
\end{remark}

\subsection{Higher KdV equations}
The examples above lead to the general problem of realizing the $(m,d+1)$-KdV evolution as the continuous limit of a $\chi$-pentagram map on curves in $\RP^d$ for any $m<d$. 
Note that such realizations (in particular, integer realizations) allow for a discretization of the KdV flow, hence the importance of and interest in this question. 
This is discussed in~\cite{Mari-Beffa:2013}, where it is conjectured that the $(m,d+1)$-flow can be realized through the intersection of an $(m-1)$-dimensional space with $m-1$ spaces of dimension $d-1$.
To approach this problem, according to Proposition \ref{prop:gamma_evolution_to_kdv},
one needs to specify the conditions which any general $\chi$ must satisfy in order for the $\e$-expansion of $\ge$ to be of the form
\begin{equation}\label{eq:higher_kdv_chi_pent}
    \ge=\Gamma+\epsilon^m c\cdot Q_m\Gamma+\O(\e^{m+1}),
\end{equation}
i.e.~$G_1=\dots=G_{m-1}=0$ and $G_m=c\cdot Q_m$ for some constant $c$.

\begin{prop}\label{prop:degrees_of_freedom_count}
    In order for the $\e$-expansion of $\ge$ to be of the form prescribed by $\eqref{eq:higher_kdv_chi_pent}$, one must impose at least $\tfrac{1}{6}(m^3-6m^2+17m-12)$ restrictions on the points in $\chi$.
\end{prop}

\begin{remark}
    This result provides a heuristic evidence against the conjecture of~\cite{Mari-Beffa:2013}. Indeed, for sufficiently large $d$, one expects that it will be impossible to realize the $(d,d+1)$-KdV equation through this way, since the number of degrees of freedom (i.e.~number of points in $\chi$) grow as~$d^2$, while the number of constraints grows as $d^3$.
    In fact, this estimate suggests that $d=9$ is already large enough for the conjecture to fail.
    Nevertheless, this proposition does not  exclude the possibility of higher KdV equations appearing as special degenerate cases.
    
    On the other hand, one would expect that if the number of degrees of freedom exceeds the number of restrictions required to make the $(m,d+1)$-KdV equation appear as the continuous limit of the $\chi$-pentagram map, for some $m$ sufficiently small, then there should exist an appropriate configuration of $\chi$ which does indeed yield this equation.
\end{remark}

\begin{proof}[Proof of Proposition~\ref{prop:degrees_of_freedom_count}]
    Each requirement $G_i=0$ corresponds to at least $\frac{1}{2}(i^2-3i+4)$ equations to be satisfied by the points in $\chi$. Indeed, demanding that $\alpha_{i,i-j}=0$ imposes one restriction when $j=0$ and none when $j=1$ (by Corollary~\ref{cor:coef_properties}). When $j\geq 2$, the coefficient $\alpha_{i,i-j}$ will be a polynomial in the $u_\ell$ functions and their derivatives. In particular, it will contain a $u_{d-1-k}^{(j-2-k)}$ term for each $0\leq k\leq j-2$, corresponding to $j-1$ restrictions. Therefore, demanding $G_i=0$ imposes at least $1+0+1+2+\dots+(i-2)=\frac{1}{2}(i^2-3i+4)$ restrictions (note that $\alpha_{i,0}$ is uniquely determined by the normalization). It then follows after summation that taking $G_1=\dots=G_{m-1}=0$ requires $\tfrac{1}{6}(m^3-6m^2+17m-12)$ restrictions.
\end{proof}

\end{document}